\documentclass[11pt,a4paper]{article}
\usepackage{amsmath}
\usepackage{amsthm}
\usepackage{amssymb}
\usepackage{mathrsfs}
\usepackage{color,eucal}

\definecolor{ddmagenta}{rgb}{0.7,0,0.9}
\definecolor{ddcyan}{rgb}{0,0.2,1.0}
\definecolor{dred}{rgb}{.8,0,0}

\newcommand{\bele}{\begin{lemm}}
\newcommand{\enle}{\end{lemm}}
\newcommand{\bedef}{\begin{defi}}
\newcommand{\bete}{\begin{teor}}
\newcommand{\eddef}{\end{defi}}
\newcommand{\ente}{\end{teor}}
\newcommand{\beos}{\begin{osse}}
\newcommand{\eddos}{\end{osse}}
\newcommand{\bepr}{\begin{prop}}
\newcommand{\empr}{\end{prop}}
\newcommand{\bepro}{\begin{prob}}
\newcommand{\empro}{\end{prob}}
\newcommand{\bede}{\begin{defin}}
\newcommand{\edde}{\end{defin}}
\newcommand{\beco}{\begin{coro}}
\newcommand{\enco}{\end{coro}}

\newcommand{\beeq}[1]{\begin{equation}
 \label{#1}}
\newcommand{\eddeq}{\end{equation}}

\newcommand{\beeqa}[1]{\begin{eqnarray}
  \label{#1}}
\newcommand{\eddeqa}{\end{eqnarray}}

\newcommand{\beal}[1]{\begin{align}
 \label{#1}}
\newcommand{\eddal}{\end{align}}

\newcommand{\bespl}[1]{\begin{split}
 \label{#1}}
\newcommand{\edspl}{\end{split}}

\newcommand{\bega}[1]{\begin{gather}
 \label{#1}}
\newcommand{\edga}{\end{gather}}

\newcommand{\beeqax}{\begin{eqnarray*}}
\newcommand{\eddeqax}{\end{eqnarray*}}

\newcommand{\tensore}{\varepsilon({\bf u})}
\newcommand{\tensoret}{\varepsilon({\mathbf{u}_t})}

\newcommand{\teta}{\vartheta}

\newcommand{\nn}{{\bf n}}
\newcommand{\uu}{{\bf u}}

\newcommand{\vv}{{\bf v}}
\newcommand{\yy}{{\bf y}}
\newcommand{\eeta}{{\mbox{\boldmath$\eta$}}}
\newcommand{\WW}{{\bf W}}

\newcommand{\eps}{\varepsilon}
\newcommand{\epsi}{\varepsilon}

\newcommand{\weak}{\rightharpoonup}

\DeclareMathOperator{\dive}{div}

\let\TeXchi\chi
\def\chi{{\setbox0 \hbox{\mathsurround0pt
$\TeXchi$}\hbox{\raise\dp0 \copy0 }}}

\newtheorem{maintheorem}{Theorem}
\newtheorem{theorem}{Theorem}[section]

\newtheorem{lemma}{Lemma}[section]
\newtheorem{proposition}[lemma]{Proposition}
\newtheorem{remark}[lemma]{Remark}%

\newtheorem{notation}[lemma]{Notation}

\renewcommand{\part}{\partial_t}

\newcommand{\nnu}{\nonumber}

\newcommand{\weaksto}{{\rightharpoonup^*}}
\newcommand{\weakto}{\rightharpoonup}

\newcommand{\pairing}[4]{ \sideset{_{#1 }}{_{ #2}}  {\mathop{\langle #3 , #4  \rangle}}}

 \def\fin{\hfill
         \trait .3 5 0
         \trait 5 .3 0
         \kern-5pt
         \trait 5 5 -4.7
         \trait 0.3 5 0
 \medskip}
 \def\trait #1 #2 #3 {\vrule width #1pt height #2pt depth #3pt}
\newcommand{\forae}{\text{for a.e.}}
\newcommand{\aein}{\text{a.e.\ in}}

\newcommand{\R}{\Bbb{R}}
\newcommand{\N}{\Bbb{N}}

\newcommand{\piecewiseConstant}[2]{\overline{#1}_{\kern-1pt#2}}

\newcommand{\underlinepiecewiseConstant}[2]{\underline{#1}_{\kern-1pt#2}}

\newcommand{\piecewiseLinear}[2]{{#1}_{\kern-1pt#2}}

\newcommand{\pwM}[2]{\widetilde{#1}_{\kern-1pt#2}}
 \def\trait #1 #2 #3 {\vrule width #1pt height #2pt depth #3pt}
\newcommand{\pwN}[2]{#1_{\kern-1pt#2}}
 \def\trait #1 #2 #3 {\vrule width #1pt height #2pt depth #3pt}

\newcommand{\dom}{\text{{\rm D}}}

\newcommand{\matrid}{\mathbf{1}}
\newcommand{\Id}{{\rm Id}}
\newcommand{\bfw}{\mathbf{W}}
\newcommand{\sig}{\sigma}

\newcommand{\hunoc}{H^1 (\Gamma_c)}
\newcommand{\ellhunoc}{\ell_{H^1 (\Gamma_c)}}

\newcommand{\ellv}{\ell}
\newcommand{\gammav}{\gamma}

\newcommand{\jvstar}{j^*}

\newcommand{\jvstarhunoc}{j^*}
\newcommand{\jvstareps}{j^*_{\mu}}
\newcommand{\jvstarhunocps}{j^*_{\mu}}
\newcommand{\gammaveps}{\gamma_\mu}
\newcommand{\gammahunoceps}{\gamma_\mu}
\newcommand{\rhoveps}{\rho_\mu}

\newcommand{\calelleps}{\mathscr{L}_\mu}
\newcommand{\calelleveps}{\mathscr{L}_\mu}

\newcommand{\pepsmu}{(\mathbf{P}_{\eps}^\mu)}
\newcommand{\pmu}{(\mathbf{P}^\mu)}
\newcommand{\dd}{\mathrm{d}}
\newcommand{\tetazero}{\teta^0}
\newcommand{\tetaessezero}{\teta_s^0}
\newcommand{\tetazep}{\teta^0_{\mu}}
\newcommand{\tetazeps}{\teta_{s,\mu}^0}
\newcommand{\wzeromu}{w_\mu^0}
\newcommand{\zzeromu}{z_\mu^0}
\newcommand{\tetame}{\teta_{\eps\mu}}
\newcommand{\tetasme}{\teta_{s,\eps\mu}}
\newcommand{\wme}{w_{\eps\mu}}
\newcommand{\zme}{z_{\eps\mu}}
\newcommand{\ume}{\uu_{\eps\mu}}
\newcommand{\chime}{\chi_{\eps\mu}}
\newcommand{\xime}{\xi_{\eps\mu}}
\newcommand{\eetame}{\eeta_{\eps\mu}}
\newcommand{\uumu}{\uu_{\mu}}
\newcommand{\chimu}{\chi_{\mu}}
\newcommand{\wmu}{w_{\mu}}
\newcommand{\zmu}{z_{\mu}}
\newcommand{\tetamu}{\teta_{\mu}}
\newcommand{\tetasmu}{\teta_{s,\mu}}
\newcommand{\ximu}{\xi_{\mu}}
\newcommand{\eetamu}{\eeta_{\mu}}
\newcommand{\tetazeroemu}{\teta_{\eps\mu}^0}
\newcommand{\tetaessezeroemu}{\teta_{s,\eps\mu}^0}

\newcommand{\OLDSTUFF}[1]{}

\newcommand{\RCOMM}[1]{{\color{red} #1}}
\renewcommand{\RCOMM}[1]{}

\begin{document}

\date{}

                                %
\title{\Large Thermal effects in adhesive contact: modelling and
analysis}
                           %
\author{
Elena Bonetti\footnote{\emph{Dipartimento di Matematica ``F.\
Casorati'', Universit\`a di Pavia. Via Ferrata 1, I--27100 Pavia,
Italy. Email: {\tt elena.bonetti\,@\,unipv.it}}}\,, \, Giovanna
Bonfanti\footnote{ \emph{Dipartimento di Matematica, Universit\`a di
 Brescia.  Via Valotti 9, I--25133 Brescia, Italy.}
 E-mail: {\tt bonfanti\,@\,ing.unibs.it}}\,,  \,
    Riccarda Rossi\footnote{
\emph{Dipartimento di Matematica, Universit\`a di
 Brescia. Via Valotti 9, I--25133 Brescia, Italy.}
 E-mail: {\tt riccarda.rossi\,@\,ing.unibs.it}}
  }

\maketitle

 \numberwithin{equation}{section}

\begin{abstract}
\noindent In this paper, we consider a contact problem with adhesion
between a viscoelastic body and a rigid support, taking thermal
effects into account. The PDE system we deal with is derived within
the modelling approach proposed by \textsc{M. Fr\'emond} and, in
particular, includes the  entropy balance equations, describing the
evolution of the temperatures  of the body and of the adhesive
material. Our main result consists in showing the existence of
global in time solutions (to a suitable variational formulation) of
the related initial and boundary value problem.
\end{abstract}

                                %
                                %
\pagestyle{myheadings} \markright{\it Bonetti, Bonfanti, Rossi /
Thermal effects in contact with adhesion}
                                %
                                %

\section{Introduction}
This paper addresses the analysis of adhesive contact between a
viscoelastic body and a rigid support, in the case when thermal
effects are included.
 Contact with adhesion is described using the modelling approach proposed by
 Fr\'emond (see \cite[Chap. 14]{fre}), which was  originally introduced for the
 isothermal case, combining the theory of damage
(see, e.g., \cite{fremond-nedjar}, \cite{bonetti-schimperna},
\cite[Chap. 12]{fre}) with the theory of unilateral contact. Indeed,
although the unilateral contact theory (which prescribes the
impenetrability condition between the bodies) does not allow for any
resistance to tension, in the adhesion phenomenon resistance to
tension is given by micro-bonds on the contact surface,  preventing
 separation. Adhesion is active if these bonds (one may think of a
``glue'' on the contact surface) are not damaged. Thus,
 the description of this phenomenon has to
take into account the state of the adhesive bonds (through a
``damage parameter'') and the microscopic movements breaking them,
as well as macroscopic deformations and displacements.

In  the recent papers \cite{bbr1} and \cite {bbr2}, we have
introduced the model and derived the corresponding initial and
boundary value problem in the isothermal case. The resulting PDE
system couples an equation for  macroscopic  deformations of the
body and a ``boundary'' equation on the contact surface, describing
the evolution of the state of the  glue by a surface damage
parameter. The system  is highly nonlinear, mainly due to the
presence of nonlinear boundary conditions and nonsmooth constraints
on the physical variables. In \cite{bbr1},   existence of a global
in time solution for a weak version of the corresponding PDE system
was proved in the case of irreversible damage dynamics for the glue.
In the subsequent contribution~\cite{bbr2},  focusing on the
reversible case, we proved well-posedness results and further
investigated the long-time behaviour of the solutions.

In this paper, we aim to  generalize   the model introduced in
\cite{bbr1} and \cite{bbr2}, including thermal effects both  on the
contact surface, and in the interior.
 We believe this  to be interesting from the modeling
perspective, because external thermal actions can in fact influence
the state of the adhesive material, see~\cite{fre}.

 In extending the
model to the non-isothermal case, we shall adopt the following
viewpoint: we shall assume that the body temperature and the glue
temperature may be different and thus governed by two distinct
entropy balance laws.

\subsection{The model and the PDE system}

\newcommand{\traccia}{{|_{\Gamma_c}}}

Let us now introduce the model and derive the corresponding initial
and boundary value problem. On a time interval $(0,T)$, we
investigate the mechanical evolution of a thermoviscoelastic body
located in a smooth bounded domain $\Omega\subset{\R}^3$, whose
boundary is
$\partial\Omega=\overline\Gamma_1\cup\overline\Gamma_2\cup\overline\Gamma_c$.
Here $\Gamma_i$, $i=1,2,c$, are open subsets in the relative
topology of $\partial\Omega$, each of them with a smooth boundary
and disjoint one from each other.  In particular, $\Gamma_c$ is the
contact surface. Hereafter we shall suppose that both $\Gamma_c$ and
$\Gamma_1$ have   positive measure. Without loss of generality , we
shall treat $\Gamma_c$ as a flat surface and identify it with a
subset of ${\R}^2$.

The thermomechanical equilibrium of the system is described by the
state variables.
 We consider the absolute temperature of the
body $\teta$ and the symmetric linearized strain tensor
$\epsilon(\uu)$ (we shall denote by $\uu$ the vector of small
displacements), defined in $\Omega\times (0,T)$. Moreover, the
variables describing the equilibrium on $\Gamma_c\times (0,T)$ are
the absolute temperature of the glue $\teta_s$, a damage parameter
$\chi$, its gradient $\nabla\chi$, and the trace
$\uu_{|_{\Gamma_c}}$ of the displacement $\uu$ on the contact
surface. The parameter $\chi$ is assumed to take values in
$[0,1]$, with $\chi=0$ for completely damaged bonds, $\chi=1$ for
undamaged bonds, and $\chi\in(0,1)$ for partially damaged bonds.

The free energy of the system is given by a volume contribution
$\Psi_\Omega$ and a surface one $\Psi_{\Gamma_c}$. It is known from
thermodynamics that the free energy is concave with respect to the
temperature. Thus,  considering   a  fairly general expression for
the purely
thermal contribution in the free energy (cf. \cite{BFR}) and 
normalizing some physical constants, we assume in $\Omega\times
(0,T)$
\begin{equation}\label{energiacorpo}
 \Psi_\Omega=-j(\teta)+p(\teta)\text{tr}(\epsilon(\uu))+
 \frac 1 2 \epsilon(\uu)K\epsilon(\uu),
\end{equation}
{where $j$ is a sufficiently regular, increasing, and convex real
function,} the function $p$ accounts for the thermal expansion
energy, and  $K=(a_{ijkh})$ denotes the elasticity tensor for a
possibly anisotropic and inhomogeneous material.  In the sequel, for
the sake of simplicity, we assume $p(\teta)=\teta$. A few comments
on the function $j$ are now in order. A possible choice for $j$,
often used in the literature, is
\begin{equation}\label{scelta}
j(\teta)=\teta\log\teta-\teta.
\end{equation}
This  enforces  the physical constraint that $\teta$ be  strictly
positive. However, in our mathematical analysis we are going to
tackle more general situations. In particular, we shall not require
any condition  on the domain of $j$. Analogously, we prescribe in
$\Gamma_c\times (0,T)$
\begin{equation}\label{energiabordo}
 \Psi_{\Gamma_c}=-j(\teta_s)+\lambda(\chi)(\teta_s-\teta_{eq})+
 I_{[0,1]}(\chi)+\sigma(\chi)+\frac 1 2|\nabla\chi|^2+
 \frac 1 2 \chi|\uu_\traccia|^2+I_-(\uu_\traccia\cdot{\bf n}),
\end{equation}
where $\teta_{eq}>0$ is a critical temperature and  $\lambda$ is a
regular (quadratic) function.  Once we
 consider  contact with adhesion as the effect of a phase
transition between the undamaged and damaged state of the adhesive
substance on the contact surface, $ \lambda'$
 formally  corresponds to the so-called latent heat in phase
transitions models and $ \teta_{eq}$  to the critical temperature
between undamaged and damaged adhesion. This  relates  to the
assumption that just by temperature devices we can damage the
micro-bonds  on the contact surface.
 Moreover, the indicator function $I_{[0,1]}$ of the
interval $[0,1]$ accounts for physical constraints on $\chi$, being
$I_{[0,1]}(\chi)=0$ if $\chi\in [0,1]$ and $I_{[0,1]}(\chi)=+\infty$
otherwise. Analogously, denoting by $I_-$ the indicator function of
the interval $(-\infty,0]$, the term $I_-(\uu_\traccia\cdot{\bf n})$
renders the impenetrability condition on the contact surface, as it
enforces that $\uu_\traccia\cdot{\bf n}\leq 0$ (${\bf n}$ is the
outward unit normal vector to $\Gamma_c$). Finally, the function
$\sigma$ is sufficiently smooth and possibly nonconvex, being
related to nonmonotone dynamics for $\chi$
 (from a physical point of view, it
corresponds to some cohesion in the material).

The free energy describes the thermomechanical equilibrium of the
system in terms of fixed state variables. Hence, we follow the
approach
 proposed by \textsc{J.J. Moreau} to prescribe the
  dissipated energy by means of a dissipation functional, the so-called pseudo-potential of
  dissipation,
which is a convex, nonnegative functional, attaining its minimum 0
when the  dissipation (described by the dissipative variables) is
zero.  The dissipative variables defined in $\Omega\times (0,T)$ are
$\nabla\teta$ and $\epsilon(\uu_t)$.
  Thus,
we define the volume part $\Phi_\Omega$ of the pseudo-potential of
dissipation by
\begin{equation}\label{pseudocorpo}
\Phi_\Omega=\frac 1 {2}|\nabla\teta|^2+\frac 1 2 \epsilon(\uu_t)
K_v\epsilon(\uu_t),
\end{equation}
where $K_v=(b_{ijkh})$ denotes the viscosity tensor for a possibly
anisotropic and inhomogeneous material. The surface part
$\Phi_{\Gamma_c}$ of the pseudo-potential of dissipation depends on
$\nabla\teta_s$, $\chi_t$, and also on the difference
$(\teta_\traccia-\teta_s)$ between the temperatures of the body and
of the glue on the contact surface, namely
\begin{equation}\label{pseudobordo}
\Phi_{\Gamma_c}=\frac 1 {2}|\nabla\teta_s|^2+\frac 1 2 |\chi_t|^2+
\frac 1 2  k(\chi) (\teta_\traccia-\teta_s)^2.
\end{equation}
Here, $k$ is a sufficiently regular function and its physical
meaning is related to the heat exchange between the body and the
adhesive material. It is of fairly natural evidence (see
also~\cite{bruneaux-benzarti08,theillout83,troung90}  that the
possibility (and the amount) of heat exchange between the body and
the contact surface depends on the fact that the adhesion is more or
less active. We let $k$ to be nonnegative (in accordance with
thermodynamical consistency ensured by the convexity of the
pseudo-potential of dissipation), increasing, and possibly vanishing
when $\chi$ attains its minimum value 0. Indeed, we may think that
if the adhesion is not active no heat exchange is allowed ($k(0)=0$)
or that a residual heat exchange is preserved even for the
completely damaged adhesive substance $(k(0)>0$).

Now, let us introduce the equations  in accordance with the laws of
thermomechanics. We consider the momentum balance (in the
quasi-static case)
\begin{equation}\label{bilanciomacro}
-\dive \Sigma={\bf f} \quad\text{in } \Omega\times (0,T),
\end{equation}
where $\Sigma$ is the stress tensor, combined with the  boundary
conditions (${\bf R}$ is the reaction on the contact surface)
\begin{equation}
\Sigma{\bf n}={\bf R}\quad\text{in }\Gamma_c\times (0,T),\quad
{\uu}={\bf 0}\quad\text{in }\Gamma_1\times (0,T),\quad\Sigma{\bf
n}={\bf g}\quad\text{in }\Gamma_2\times (0,T),
\end{equation}
${\bf f}$ being a volume force and ${\bf g}$ a traction. The thermal
balance   is given by the following entropy equation
\begin{equation}\label{eqentropiacorpo} s_t+\dive{\bf
Q}=h\quad\text{in }\Omega\times (0,T),
\end{equation}
$s$ denoting the internal entropy, ${\bf Q}$ the entropy flux, and
$h$ and external entropy source.
 Indeed, equation \eqref{eqentropiacorpo} can
be obtained rescaling the first law of thermodynamics (dividing the
internal energy balance  by the absolute temperature), under the
small perturbation assumption. We refer to \cite{bcf}, \cite{BFR} and
 \cite{bcfg1}, \cite{bcfg2}, \cite{bcfgIII} for details on this modelling
approach and related analytical results. Moreover, we supplement
\eqref{eqentropiacorpo} with the following boundary conditions
\begin{equation}\label{flussoentropia}
{\bf Q}\cdot{\bf n}={F}\,\text{ on }\Gamma_c\times (0,T),\quad{\bf
Q}\cdot{\bf n}=0\,\text{ on }\partial\Omega\setminus\Gamma_c\times
(0,T),
\end{equation}
where ${F}$ denotes the entropy flux through $\Gamma_c$. Hence, on
the contact surface, we introduce a balance equation for
 the microscopic movements,
 that is
\begin{equation}\label{bilanciomicro}
 B-\dive {\bf H}=0\,\text{ in }\Gamma_c\times (0,T),\quad{\bf H}\cdot{\bf n}_s=0
\,\text{ on }\partial\Gamma_c\times (0,T),
\end{equation}
$B$, ${\bf H}$ representing interior forces, responsible
for the damage of adhesive bonds between the body and the support,
and ${\bf n}_s$ the outward unit normal vector to
$\partial\Gamma_c$. Then, the entropy equation on the contact
surface is given by
\begin{equation}\label{entropiabordo}
 \partial_t s_s+\dive{\bf Q}_s={F}\,\text{ in }\Gamma_c\times (0,T),
 \quad{\bf Q}_s\cdot{\bf n}_s=0\,\text{ on }\partial\Gamma_c\times (0,T).
\end{equation}
Here, $s_s$ is the entropy on the contact surface, ${\bf Q}_s$ the
surface entropy flux, and the term $F$, given by the flux through
the boundary $\Gamma_c$ (cf. \eqref{flussoentropia}), represents a
surface entropy source.

Constitutive relations for $\Sigma, {\bf R}, s, {\bf Q}, F, B, {\bf
H}, s_s, {\bf Q}_s$ are given in terms of the free energies and the
pseudo-potentials of dissipation. More precisely, we have {
\begin{align}
 &s=-\frac{\partial\Psi_\Omega}{\partial\teta}=\ell(\teta)-\dive{\uu},\\
&s_s=-\frac{\partial\Psi_{\Gamma_c}}{\partial\teta_s}=\ell(\teta_s)-\lambda(\chi),
\end{align}
where $\ell$ is the derivative} of the convex function $j$. In the
physical case $j(x)=x\log x-x$ (cf. \eqref{scelta}), we have
\begin{equation}\label{fisico}
\ell(x)=\log x.
\end{equation}
{In particular,  $\ell$   in \eqref{fisico} in fact yields an
internal positivity constraint on the system temperatures $\teta$
and $\teta_s$.}   Furthermore,
\begin{align}
 &{\bf Q}=-\frac{\partial{\Phi_\Omega}}{\partial\nabla\teta}=-\nabla\teta,\\
 &{\bf Q}_s=-\frac{\partial{\Phi_{\Gamma_c}}}{\partial\nabla\teta_s}=-\nabla\teta_s,\\
 & F=\frac{\partial{\Phi_{\Gamma_c}}}{\partial(\teta_\traccia-\teta_s)}=k(\chi)(\teta_\traccia-\teta_s).
\end{align}
The constitutive relation for the stress tensor $\Sigma$ accounts
for dissipative (viscous) dynamics for deformations, in that we have
\begin{equation}
\Sigma=\frac{\partial\Psi_\Omega}{\partial\epsilon(\uu)}+
\frac{\partial\Phi_\Omega}{\partial\epsilon(\uu_t)}=K\epsilon(\uu) +
K_v\epsilon(\uu_t)+\teta {\bf 1},
\end{equation}
(${\bf 1}$ denotes the identity matrix), while the reaction ${\bf
R}$ is given by
\begin{equation}
{\bf
R}=-\frac{\partial\Psi_{\Gamma_c}}{\partial\uu_\traccia}=-\chi\uu_\traccia
-\partial I_-(\uu_\traccia\cdot {\bf n}){\bf n}.
\end{equation}
 We further prescribe $B$
\begin{align}
&B=\frac{\partial\Psi_{\Gamma_c}}{\partial\chi}+\frac{\partial\Phi_{\Gamma_c}}{\partial\chi_t}\\\nonumber
&=\lambda'(\chi)(\teta_s-\teta_{eq})+\partial
I_{[0,1]}(\chi)+\sigma'(\chi)+\frac 1 2 |\uu_\traccia|^2+\chi_t,
\end{align}
($\partial I_-$ and $\partial I_{[0,1]}$ standing for the
subdifferentials of $I_-$ and $I_{[0,1]}$, respectively), and let
${\bf H}$ be
\begin{equation}
{\bf
H}=\frac{\partial\Psi_{\Gamma_c}}{\partial\nabla\chi}=\nabla\chi.
\end{equation}
\begin{remark}
\upshape \label{rem:elena}
 Let us point out that the evolution of the system
 is characterized by dissipation due to choice of dissipative
  potentials $\Phi_\Omega$ and $\Phi_{\Gamma_c}$ (cf.
   \eqref{pseudocorpo} and \eqref{pseudobordo})
    and the balance laws of thermodynamics (see in particular
     \eqref{eqentropiacorpo} and \eqref{entropiabordo},
       \eqref{bilanciomacro}, and \eqref{bilanciomicro}).
       \end{remark}

Combining the previous constitutive relations with the balance laws,
 we obtain
 the following  boundary value problem
\begin{align}
\label{eqteta} &\partial_t (\ell(\teta)) - \dive (\uu_t) -\Delta
\teta= h \qquad \text{in $\Omega \times (0,T)$,}
\\
\label{condteta} &\partial_n \teta= \begin{cases} 0 & \text{in
$(\partial \Omega \setminus \Gamma_c) \times (0,T)$,}
\\
-k(\chi) (\teta-\teta_s) & \text{in $ \Gamma_c \times (0,T)$,}
\end{cases}
\\
\label{eqtetas} &\partial_t (\ell(\teta_s)) - \partial_t
(\lambda(\chi))  -\Delta \teta_s= k(\chi) (\teta-\teta_s) \qquad
\text{in $\Gamma_c \times (0,T)$,}
\\
\label{condtetas} &\partial_n \teta_s =0 \qquad \text{in $\partial
\Gamma_c \times (0,T)$,}
\\
\label{eqI}
&-\hbox{div }(K\tensore+ K_v\tensoret+ \teta \matrid)={\bf f}\qquad\hbox{in }\Omega\times (0,T),\\
&{\bf u}={\bf 0}\quad\hbox{in }\Gamma_{1}\times (0,T),
\quad(K\tensore+K_v\tensoret + \teta \matrid){\bf n}={\bf g}\qquad\hbox{in }\Gamma_{2}\times (0,T),\label{condIi}\\
& (K\tensore+ K_v\tensoret + \teta \matrid){\bf n}+\chi{\bf u}+
\partial I_-(\uu\cdot {\bf n}){\bf n}
\ni{\bf 0}\qquad\hbox{in }\Gamma_{c} \times (0,T),\label{condIii}\\
&\chi_t -\Delta\chi+\partial I_{[0,1]}(\chi) +
\sigma'(\chi)-\lambda'(\chi) \teta_{\rm{eq}}
 \ni -\lambda'(\chi)\teta_s-\frac 1 2\vert{\bf u}\vert^2 \qquad\hbox{in }\Gamma_c
\times (0,T),\label{eqII}
\\
&\partial_{\nn} \chi=0 \qquad \text{in $\partial \Gamma_c \times
(0,T)\,$} \label{bord-chi}
\end{align}
(here  and in what follows, we shall omit for simplicity the index
$v_\traccia$ to denote the trace  on $\Gamma_c$ of a function $v$,
defined in $\Omega$). 
{
\begin{remark}\label{jgen}
\upshape As we have already pointed out, our analysis actually
accounts for a
 form of the thermal contribution in the free energies more general
than \eqref{scelta} (cf. \eqref{energiacorpo} and
\eqref{energiabordo}). Thus, our results may apply to several
physical situations with different thermal behaviour. We note that
we can also handle the case when the specific heat  $c_V$, given by
the thermodynamic relation $ c_V=\displaystyle {-\teta \frac
{\partial^2 \Psi_\Omega}{\partial\teta^2}}$
(cf.~\cite[Rem.~2.1]{BFR}), is not  constant. Indeed, it is known
from physics that the specific heat may be depending on the
temperature, e.g., $c_V (\teta)=\teta^\gamma$, $\gamma>0$. A choice
of $j$ corresponding to the latter expression of $c_V$ is
$$
j(\teta)=\frac{\teta^{\gamma+1}}{\gamma(\gamma+1)}\,,
$$
which is covered by our analysis. In particular, letting $\gamma=1$
(and hence $j(\teta)=\frac 1 2\teta^2$), we get $\ell(\teta)=\teta$,
so that \eqref{eqteta} and \eqref{eqtetas} reduce to  Caginalp-type
heat equations (cf. with~\cite{caginalp86}).
\end{remark}

\subsection{Related literature and our own results}
 For a review of the theory of contact problems  see, e.g., the
monographs \cite{shillor-sofonea-telega04, eck05,
sofonea-han-shillor},  and the references therein.
 We also refer to~\cite{bbr1, bbr2} for some partial survey of
the literature on (isothermal) models of adhesive contact.  In this
connection, we  mention~\cite{Point}, where a unilateral contact
model (also derived within Fr\'emond's approach)  is considered.
Therein,
 the adhesive properties are described by the
condition $\chi\uu=0$ on the contact surface. The author proves  an
existence theorem for the related PDE system and also develops some
numerical investigations. The paper  \cite{Raous} focuses on a model
combining  unilateral contact with adhesion and friction as well:
 under a smallness condition
on the friction coefficient, an existence result is proved and
various numerical schemes are proposed.

As for  the literature on  contact models including thermal effects
in the three-dimensional framework,
 besides the contributions  mentioned in~\cite{eck05}
  we may  recall \cite{Shi-Shillor},
where a frictionless contact problem between a thermoelastic body
and a rigid foundation is modelled by a parabolic equation for
$\teta$, coupled with an elliptic equation for $\uu$, with mixed
boundary conditions. An existence result is proved, provided  the
coefficient of thermal expansion is sufficiently small. Among
dynamic models, in which  contact is rendered by means of a normal
compliance condition, we quote \cite{Figue-Trabucho}, which deals
with a wide class of frictional contact problems in thermoelasticity
and thermoviscoelasticity. Moreover, we mention \cite{Andrews},
where a frictional contact problem involving a thermoelastic body
undergoing wear on the contact surface is investigated. A
well-posedness result is proved for a system coupling a parabolic
equation for the temperature, a variational inequality for the
displacement, and a first order equation for the wear function,
supplemented with nonlinear boundary conditions.
\parskip2mm

Nevertheless, as far as we know, no results are available, in the
literature, on unilateral contact models which take into account
both adhesive properties and thermal effects. Indeed, one of the
main novelties of the present contribution is that we consider heat
generation effects in the adhesive contact phenomenon, too. That is
to say, we allow for  the body and the adhesive material to have
different temperatures, whose evolution is mainly ruled by the heat
exchange throughout the contact surface. More precisely, the entropy
flux $F$ through $\Gamma_c$ (occurring in \eqref{flussoentropia})
plays the role of a source of entropy in \eqref{entropiabordo}. From
an analytical point of view,  this results in a nonlinear coupling
between \eqref{eqteta}--\eqref{condteta} and \eqref{eqtetas} and
gives raise  to some   technical difficulties. A further peculiarity
of our work consists in assuming entropy balance laws (in place of
the more usual internal energy balance), for describing the
evolution of the body and of the glue temperatures. This brings to
strong and possibly singular nonlinearities in \eqref{eqteta} and
\eqref{eqtetas} (see \eqref{fisico}).
 An advantage of this choice is  that,  assuming that the domain $D(j)\subseteq(0,+\infty)$ (as
 in the case of
the classical choice \eqref{scelta}) once the problem is solved in a suitable
 sense,
 the positivity of the
temperatures is deduced. On this fact the thermodynamical
consistency of the model relies, see also~\cite{bcf, bcfg1, bcfg2,
BFR}. This is of particular interest in the present case, since the
low spatial regularity of the solution components $\teta$ and
$\teta_s$, along with the nonlinear boundary
condition~\eqref{condteta}, prevents us from using any maximum
principle.
\parskip0mm

In fact,  we shall study the Cauchy problem  for a {\em generalized
version} of  system \eqref{eqteta}--\eqref{bord-chi}, see
Problem~$(\mathbf{P})$ in Sec.~\ref{s:2.2}. Namely, we replace the
subdifferential operators in~\eqref{condIii}--\eqref{eqII}, by
general maximal monotone operators (possibly  rendering physical
constraints on the variables $\chi$ and $\uu$). {Further,
 we generalize the choice of  the nonlinearity $\ell$, allowing for a   maximal monotone
operator. In fact, the only restriction we impose on $\ell$ is that
the resulting internal energy of the system be coercive, cf.
with~\eqref{e:coerc-intro} below. This
 is
reasonable from a physical point of view and still enables us to
include several choices of $\ell$ in our analysis
 (in particular,  \eqref{scelta},  as well as the examples of Remark~\ref{jgen}). The idea is
 that,
once the internal energy of the system is bounded, the absolute
temperature  is bounded, too. Let us focus on the volume temperature
$\teta$.  As  known from thermodynamics, the internal energy
(depending on the entropy $s$) can be introduced as the convex
conjugate function (with respect to the variable $\teta$), of the
negative
 of the free energy (which is convex w.r.t. the temperature as
the function $j$ is convex, see \eqref{energiacorpo} and
\eqref{energiabordo}). Namely, in the case of the volume free energy
the related internal energy is
$$
e(s,\cdot)=(-\Psi_\Omega)^*(s,\cdot)=\sup_{\teta}(s\teta+\Psi_\Omega(\teta,\cdot)).
$$
Thus, our coercivity condition may be expressed  in terms of the
conjugate of $j$ by
\begin{equation}
\label{e:coerc-intro} \exists\, C_1, C_2>0: \quad  j^*(y)\geq
C_1|x|-C_2 \quad \text{if }y=\ell(x)
\end{equation}
(see subsequent \eqref{hyp:gamma-2-a} and Remark \ref{rem:jstar}).}

The main difficulties attached to the analysis of the  PDE system
\eqref{eqteta}--\eqref{bord-chi}
 are related to  the
singular character of the entropy equations~\eqref{eqteta}
and~\eqref{eqtetas}, to the nonlinear coupling between the latter
equations,  as well as between \eqref{condIii} and \eqref{eqII},
and, last but not least, to the presence of {\em general}
 multivalued
operators in all of the latter equations. In particular, it seems to
us that dealing with a general maximal monotone operator
in~\eqref{eqteta} and~\eqref{eqtetas} brings about some technical
difficulties, particularly in connection  with the third type
boundary condition~\eqref{condteta} for $\teta$ on $\Gamma_c$.

 All of these peculiarities will be carefully handled in the
 proof of our main result, Theorem 1 (see Sec.~\ref{s:2.2}), stating
 the existence of {\em global in time} solutions to the Cauchy
 problem for (the generalized version of) system
 \eqref{eqteta}--\eqref{bord-chi}. We sketch below the main steps of
 our procedure, based on a suitable approximation of
 Problem~$(\mathbf{P})$,   and on the derivation of suitable a priori estimates, which
 enable us to pass to the limit in the approximation. Such estimates are intrinsically related
to the dissipative character of the system, highlighted in
Remark~\ref{rem:elena}. In fact, in the paper \cite{bbr4} we take
advantage of  the dissipative character of the system to perform its
long-time analysis, showing that in the limit we reach a stationary
equilibrium in which dissipation vanishes.

 Note, however, that  uniqueness  is still an open issue,  at
least in the functional framework of our existence theorem. Without
going into details, we may point out that the major
 obstacle  is due to the singular character of
equations \eqref{eqteta} and \eqref{eqtetas}. In particular, the
boundary condition~\eqref{condteta} makes it harder
 to prove contraction
estimates leading to uniqueness. 
We refer to Remark~\ref{rem:explanation} for additional observations
on this point.
 Actually, uniqueness holds in
the (more regular) framework of the approximate problem, see
Section~\ref{s:3.5}.

\paragraph{Plan of the paper.} In Section
2, we  enlist our assumptions on the problem data, present
 the variational formulation of the Cauchy problem for (a
generalized version of) system \eqref{eqteta}--\eqref{bord-chi},
and
 state of our main result.
 In Section 3,
 we set up the approximation of   Pro\-blem~$(\mathbf{P})$,
suitably regularizing the maximal monotone operators in equations
\eqref{eqteta} and \eqref{eqtetas} and therein inserting (vanishing)
viscosity terms. Hence, we prove  a well-posedness result for the
approximate problem. We combine a Schauder fixed point technique for
local existence with a prolongation argument, based on global in
time {\it a priori} estimates, while uniqueness follows from   contraction
estimates. Next, in
 Section 4 we pass to  the limit in the approximate problem
  by compactness and monotonicity tools, and show that the
  approximate solutions converge to a solution of
  Problem~$(\mathbf{P})$.
 Finally, in  the Appendix  we prove some auxiliary technical results.

\section{Main result}
\label{s:2}
\subsection{Setup and preliminary results}
\label{s:2.1}
\begin{notation}
\upshape \label{n:1}
 Throughout the paper,
   given a Banach space $X$, we   denote by
 $_{X'}\langle\cdot,\cdot\rangle_X$ the duality pairing
between  $X'$ and $X$ itself, and by $\Vert\cdot\Vert_X$  both the
norm in  $X$ and in any power of $X$;  $C_{\text{w}}^0([0,T];X)$ is
the space of weakly continuous $X$-valued functions on $[0,T]$.
Whenever $X= Y_1 \times \ldots \times Y_N$, we denote by $\pi_i$,
$i=1, \ldots, N$ the projection on the $i$-th component.
\end{notation}
%
\paragraph{Young inequalities.}
We recall the Young inequality for convolutions, namely
\begin{equation}
\label{young}
\begin{gathered}
\forall\, p, \, q, \, r \in [1, \infty]
 \ \text{s.t} \ \  \frac1r=\frac1p+\frac1q-1
\quad \forall\, a \in L^p (0,T) \ b \in L^q (0,T;X) \quad \text{we
have}
\\
a * b \in L^r (0,T;X) \quad \text{and} \ \ \| a * b \|_{L^r (0,T;X)}
\leq \| a \|_{L^p (0,T)} \| b \|_{L^q (0,T;X)}\,,
\end{gathered}
\end{equation}
 and the
Young inequality
\begin{equation}
\label{young-p-q} \forall\, \delta >0 \ \ \exists \, C_\delta>0 \, :
\quad  \forall\, p, \, q \in (1,\infty) \ \text{with} \
\frac1p+\frac1q=1 \ \ \  ab \leq \delta a^p + C_\delta b^q \qquad
\text{for all $\, a, \, b \in \R$}\,.
\end{equation}
\paragraph{Functional setup.} Henceforth,  we shall suppose that
  $\Omega$ is a
bounded smooth set of $\R^3$, such that $\Gamma_c$ is a smooth
bounded domain of $\R^2$, and use the notation
$$
\begin{gathered}
H:=L^2(\Omega), \quad V:=H^1(\Omega), \quad \text{and}
\\
\bfw:= \left\{\vv \in V^3 \, : \ \vv={\bf 0} \, \hbox{ a.e. on
}\Gamma_1\right\}\,,
\end{gathered}
$$
the latter space endowed with the norm induced by $V$.
 {We shall  work with    the standard Riesz operator
\begin{equation}
\label{e:riesz}
 \mathcal{R}: V \to V' \ \ \text{given by}  \ \
 \pairing{V'}{V}{\mathcal{R}(u)}{v}:= \int_{\Omega} uv + \int_{\Omega} \nabla u \nabla
 v \quad \text{for all $u,\, v \in V$,}
\end{equation}
  and denote by   $\mathcal{R}_{\Gamma_c}$ the analogously defined  Riesz operator mapping $H^1
(\Gamma_c)$ into $(H^1 (\Gamma_c))'$. Further, we shall} extensively
use that
\begin{align}
& \label{continuous-embedding} V \subset L^p (\Gamma_c) \ \
\text{with a continuous (compact) embedding for} \ \ 1 \leq p \leq 4
\ \ \text{($1 \leq p <4$, resp.)},
\\
& \label{continuous-embedding-2} H^1 (\Gamma_c) \subset L^p
(\Gamma_c)\ \  \text{with a compact embedding for}\ \  1 \leq p<
\infty.
\end{align}
For simplicity, we denote by $\int_{\Gamma_c}uv$
($\int_{\Gamma_2}uv$, resp.) the duality pairing $
\pairing{(H^{-1/2}(\Gamma_c))^3}{(H^{1/2}(\Gamma_c))^3}{u}{v}$
between $(H^{-1/2}(\Gamma_c))^3$ and $(H^{1/2}(\Gamma_c))^3$
(between $(H^{-1/2}(\Gamma_2))^3$ and $(H^{1/2}(\Gamma_2))^3$,
resp.).  Finally, given a subset $\mathcal{O} \subset \R^N$,
$N=1,2,3,$ we shall denote by $|\mathcal{O}|$ its Lebesgue measure.
\paragraph{Preliminaries of  viscoelasticity  theory.}
We now introduce the standard bilinear forms of linear
viscoelasticity
 which allow us  to give a
variational formulation of equation~\eqref{eqI}. Dealing with an
anisotropic and inhomogeneous material,
 we assume that the fourth-order tensors $K=(a_{ijkh})$ and  $K_v=(b_{ijkh})$,
denoting the elasticity and the viscosity tensor, respectively,
satisfy the classical symmetry and ellipticity conditions
$$
\begin{aligned}
& a_{ijkh}=a_{jikh}=a_{khij}\,,\quad  \quad
b_{ijkh}=b_{jikh}=b_{khij}\,,
  \quad   i,j,k,h=1,2,3
\\
&  \exists \, \alpha_0>0 \,:  \qquad a_{ijkh} \xi_{ij}\xi_{kh}\geq
\alpha_0\xi_{ij}\xi_{ij} \quad   \forall\, \xi_{ij}\colon \xi_{ij}=
\xi_{ji}\,,\quad i,j=1,2,3\,,
\\
& \exists \, \beta_0>0 \,:  \qquad b_{ijkh} \xi_{ij}\xi_{kh}\geq
\beta_0\xi_{ij}\xi_{ij} \quad   \forall\, \xi_{ij}\colon \xi_{ij}=
\xi_{ji}\,,\quad i,j=1,2,3\,,
\end{aligned}
$$
where the usual summation convention is used.
Moreover, we require
$$
a_{ijkh}, b_{ijkh} \in L^{\infty}(\Omega)\,,
\quad  i,j,k,h=1,2,3.
$$
By the previous assumptions on the elasticity and viscosity coefficients,
the following bilinear forms
$a, b : \bfw \times \bfw \to \R$,   defined~by
$$
\begin{aligned}
a({\bf u},{\bf v}):=\int_{\Omega} a_{ijkh} \varepsilon_{kh}({\bf u}) \varepsilon_{ij}({\bf v})
\quad  \forall \uu, \vv \in \bfw,
\\
b({\bf u},{\bf v}):=\int_{\Omega} b_{ijkh} \varepsilon_{kh}({\bf u}) \varepsilon_{ij}({\bf v})
\quad  \forall \uu, \vv \in \bfw
\end{aligned}
$$
turn out to be continuous and symmetric. In particular, we have

\begin{equation}
\label{continuity} \exists \, M >0: \ |a(\uu, \vv)| +   |b(\uu,
\vv)| \leq M \| \uu\|_\bfw \| \vv\|_\bfw \quad \forall \uu, \vv \in
\bfw.
\end{equation}
Moreover, since $\Gamma_1$ has positive measure,
 by Korn's inequality we deduce that $a(\cdot,\cdot)$ and
$b(\cdot,\cdot)$ are $\bfw$-elliptic, i.e., there exist $C_{a},
C_{b}>0 $ such that
\begin{align}
& a({\bf u},{\bf u})\geq C_a\Vert{\bf u}\Vert^2_\bfw \qquad \forall
\, \uu\in \bfw, \,\label{korn_a}
\\
& b({\bf u},{\bf u})\geq C_b\Vert{\bf u}\Vert^2_\bfw \qquad \forall
\, \uu\in \bfw. \label{korn_b}
\end{align}
Relying on  Green's formula (see, e.g., \cite{duvaut-lions}), the
variational formulation of \eqref{eqI} (cf.  \eqref{eqIa} below) can
be derived by a standard argument.
\subsection{Statement of the assumptions}
As we mentioned in the introduction, we shall address the analysis
of {a {\em generalized version} of
system~\eqref{eqteta}--\eqref{bord-chi}, in which the operators
occurring in~\eqref{condIi} and~\eqref{eqII} are replaced by
ge\-ne\-ral maximal monotone operators.  We now enlist our
assumptions on the involved} nonlinearities and on the problem data.

We consider
\begin{equation}
\label{hyp:gamma-1} \tag{2.H1} \text{ a proper, convex, and l.s.c.
function} \
 j:\R \rightarrow(-\infty,+\infty]\,,
\end{equation}
and
\begin{equation}
\label{def-log} \text{its subdifferential in the sense of convex
analysis} \ \ell=\partial j: \R \rightarrow 2^\R\,.
\end{equation}
The crucial assumption on $j$ is that the following
\emph{coercivity} condition holds:
\begin{equation}
\label{hyp:gamma-2-a} \tag{$2.\text{H}2$}
\begin{aligned}
 \exists \, C_1, \ C_2   >0  \quad
\forall\, x \in \dom(\ell), \ y \in \ellv(x) \
  : \quad
  yx - j(x) \geq C_1 |x| -C_2\,.
  \end{aligned}
\end{equation}
\begin{remark}
\label{rem:jstar} \upshape Note that \eqref{hyp:gamma-2-a} can be
rephrased as
\begin{equation}
 \label{e:coercivity-puntuale}
\exists \, C_1, \ C_2   >0 \quad \forall\, x \in \dom(\ell), \ y \in
\ell(x) \
  : \quad
 j^*(y) \geq C_1 |x| -C_2\,.
\end{equation}
{In particular, \eqref{e:coercivity-puntuale} is  fulfilled when
$j(x)=x(\log(x) -1)$ for all $x \in (0,+\infty)$  and hence $\ell$
is the logarithmic nonlinearity, i.e.  $
 \ell(x)=\log (x)
 $ for all
$x \in (0,+\infty)$. In  this case,}  simple computations show that
$j^* (y)=e^y= \ell^{-1}(y)$ for all $y \in \R$,
whence~\eqref{e:coercivity-puntuale}.
\end{remark}
 \noindent Henceforth, we shall denote by
\[
\text{
 $\gamma$ the inverse of
the operator $\ell$ (recall that $\gamma=\partial j^*$)}
\]
and, with a slight abuse of notation, we shall call $\ell$ as well
the realization of \eqref{def-log} as a maximal monotone operator
$\ell: L^2 (0,T;H) \to 2^{L^2 (0,T;H)}$ (as a  maximal monotone
operator $\ell: L^2 (0,T;L^2 (\Gamma_c)) \to 2^{L^2 (0,T;L^2
(\Gamma_c))}$, respectively).  

\paragraph{Assumptions on the other problem nonlinearities.}
 Further, we let
\begin{equation}
\label{hyp:alpha} \tag{2.H3}
\begin{gathered}
 \widehat \alpha : (H^{1/2}(\Gamma_c))^3  \rightarrow [0, + \infty] \, \text{be a  proper,
 convex and  l.s.c. functional,}
\\
\text{with $ \widehat\alpha(\mathbf{0})=0=\min \widehat\alpha $,
and}
\\
\text{we set} \  \alpha:= \partial \widehat\alpha: \,
 (H^{1/2}(\Gamma_c))^3\to 2^{(H^{-1/2}(\Gamma_c))^3}.
\end{gathered}
\end{equation}
Indeed, $\alpha$ shall generalize the subdifferential operator
appearing in~\eqref{condIii}.
\begin{remark}
\label{suit-alpha} \upshape In fact,
 condition~\eqref{condIii} may be rendered rigorously (see~\cite[Sec.~II]{bbr1} for details) by
introducing the set
$$
{\cal X}_- :=\{\vv\in (H^{1/2}(\Gamma_c))^3: \vv \cdot \nn\leq 0
\hbox{ a.e. in }\Gamma_c \},
$$
with indicator function   $I_{{\cal X}_-}$. Then we consider the
maximal monotone operator
  $\partial I_{{\cal X}_-}:(H^{1/2}(\Gamma_c))^3\to2^{(H^{-1/2}(\Gamma_c))^3
  }$, given by
  $$
\begin{aligned} &{\bf \eeta}\in (H^{-1/2}(\Gamma_c))^3 \hbox{ belongs to }
\partial I_{{\cal X}_-}(\yy) \hbox{ if and only if }\nnu\\
&\yy \in {\cal X}_-,\quad\int_{\Gamma_c} {\bf \eeta}\cdot{(\vv-\yy)}
\leq0\quad\forall \vv\in {\cal X}_-.
\end{aligned}
$$
\end{remark}
%
\noindent In the same way, in equation~\eqref{eqII} we shall
consider
\begin{equation}
  \label{A5}\tag{2.H4}
  \text{a maximal monotone
operator}  \ \ \beta: \R \to 2^{\R}, \ \ \text{with  domain
$\dom(\beta)\subseteq[0,+\infty)$.}
\end{equation}
It is well known that there exists a proper, l.s.c. and convex
function $\widehat\beta : \overline{\dom (\beta)} \to
(-\infty,+\infty]$ such that
$\beta=\partial\widehat\beta$. 
Concerning the functions  $\lambda$, $k$, and $\sig'$,
 we assume that
\begin{equation}
\label{hyp-sig} \tag{2.H5} \sig' \, : \R \to \R \ \ \text{is
Lipschitz continuous, with Lipschitz constant $L_{\sig}$,}
\end{equation}
\begin{equation}
\label{hyp-k} \tag{2.H6}  k \, : \R \to [0,+\infty)\ \  \text{is
Lipschitz continuous , with Lipschitz constant $L_{k}$,}
\end{equation}
\begin{equation}
\label{hyp-lambda} \tag{2.H7}
  \lambda \in {\rm C}^{1,1} (\R)\,,
\end{equation}
i.e., $\lambda $ has a Lipschitz continuous derivative. As a
consequence,
\begin{align}
& \label{at-most-quadratic}
 \exists\, C_3>0 \ \forall\, x,y \in \R\ :
\ \  |\lambda(x) -\lambda(y)| \leq C_3 (|x|+|y|+1)|x-y|\,,
\\
& \label{e:allafineserve}
 \exists\, C_3'>0 \ \forall\, x \in \R\ :
\ \  |\lambda'(x)| \leq C_3' (|x|+1)\,.
\end{align}
\paragraph{Assumptions on the problem data.} \
We assume that
\begin{equation}
\label{hypo-h} \tag{2.H8}  h \in L^2 (0,T;V')\cap L^1 (0,T;H) \,,
\end{equation}
\begin{equation}
\label{hypo-f} \tag{2.H9} \mathbf{f} \in L^2 (0,T;H^3)\,,
\end{equation}
\begin{equation}
\label{hypo-g} \tag{2.H10} \mathbf{g} \in L^2 (0,T;
(H^{-1/2}(\Gamma_2))^3)\,.
\end{equation}
 Then, we remark that the function
$\mathbf{F}:(0,T) \to \bfw'$ defined  by
$$
 \pairing{\bfw'}{\bfw}{\mathbf{F}(t)}{\vv}:\int_{\Omega}\mathbf{f}(t) \cdot {\vv}
 +\int_{\Gamma_{2}} \mathbf{g}(t) \cdot {\vv}
\quad \forall\, {\vv} \in W \quad \forae \, t \in (0,T),
$$
fulfils
\begin{equation}
\label{effegrande} \mathbf{F} \in L^2(0,T;\bfw').
\end{equation}
Finally, we require that the initial data fulfil
\begin{align}
& \label{cond-teta-zero} w_0 \in H \ \  \text{with} \ \ j^*(w_0) \in
L^1 (\Omega)\,,
\\
& \label{cond-teta-esse-zero} z_0 \in L^2 (\Gamma_c) \ \ \text{and}
\ \ j^*(z_0) \in L^1 (\Gamma_c)\,,
\\
& \label{cond-uu-zero} {\bf u}_0 \in \bfw \ \ \text{and} \ \  \uu_0
\in \dom (\widehat{\alpha})\,,
\\
& \label{cond-chi-zero} \chi_0 \in \hunoc\,,  \ \
\widehat{\beta}(\chi_0)\in
 L^1(\Gamma_c)\,.
\end{align}
\subsection{Variational formulation and statement of the main result}
\label{s:2.2} We  now  state the variational formulation of the
initial-boundary value problem for a generalized version of
system~\eqref{eqteta}--\eqref{bord-chi}, featuring the
nonlinearities introduced above.
\paragraph{\bf Problem $(\mathbf{P})$.} \, Given a  quadruple of initial data $(w_0,
z_0, \uu_0, \chi_0)$ complying
with~\eqref{cond-teta-zero}--\eqref{cond-chi-zero}, find  $(  \teta,
w, \teta_s,  z, \uu,\chi,\eeta,  \xi)$ such that
\begin{align}
& \teta \in L^2 (0,T; V) \cap L^\infty (0, T;L^1 (\Omega))\,, \label{reg-teta}\\
& w \in L^\infty (0,T;H) \cap H^1 (0,T;V')\,, \label{reg-w}
\\
& j^*(w) \in L^\infty (0,T; L^1(\Omega))\,, \label{reg-w-ult}
\\
& \teta_s \in  L^2 (0,T; \hunoc) \cap L^\infty (0, T;L^1 (\Gamma_c))\,, \label{reg-teta-s}\\
&
  \label{reg-zeta}
  z \in L^\infty (0,T;L^2 (\Gamma_c)) \cap H^1
  (0,T;H^1(\Gamma_c)')\,,
\\
& \jvstarhunoc(z) \in L^\infty (0,T;L^1(\Gamma_c))\,,
\label{reg-z-ult}
  \\
&{\bf u}\in H^1(0,T;\bfw)\,,\label{reguI}\\
&  \eeta\in L^2(0,T; (H^{-1/2}(\Gamma_c))^3)\,, \label{etareg}
\\
&\chi \in L^{2}(0,T;H^2 (\Gamma_c))  \cap L^{\infty}(0,T;H^1
(\Gamma_c)) \cap H^1 (0,T;L^2 (\Gamma_c))\,, \label{regchiI}
\\
& \xi\in L^2(0,T; L^{2}(\Gamma_c))\,, \label{xireg}
\end{align}
 fulfilling  the initial
 conditions
\begin{align}
& \label{iniw} w(0)=w_0 \quad \aein \ \Omega\,,
\\
& \label{iniz} z(0)=z_0 \quad \aein \ \Gamma_c\,,
\\
& \label{inichi} \chi(0)=\chi_0 \quad {\aein \ \Gamma_c}\,,
\\
& \label{iniu} \uu(0)=\uu_0 \quad {\aein \ \Omega}\,,
\end{align}
  and
 \begin{align}
&
 \label{teta-weak}
\begin{aligned}
\pairing{V'}{V}{w_t}{v} -\int_{\Omega} \dive(\uu_t) \, v
+\int_{\Omega} \nabla \teta \, \nabla v + \int_{\Gamma_c} k(\chi) &
(\teta-\teta_s) v \\ & =  \pairing{V'}{V}{h}{v} \quad \forall\, v
\in V \ \hbox{ a.e. in }\, (0,T)\,,
\end{aligned}
\\
& \label{inclu-w} {w(x,t) \in \ell(\teta(x,t)) \quad \forae\, (x,t)
\in \Omega \times (0,T)\,,}
\\
& \label{teta-s-weak}
\begin{aligned}
\pairing{{\hunoc}'}{\hunoc}{z_t}{v} -\int_{\Gamma_c}
\partial_t \lambda(\chi) \, v   & +\int_{\Gamma_c} \nabla \teta_s  \, \nabla
v \\ &= \int_{\Gamma_c} k(\chi) (\teta-\teta_s) v  \quad \forall\, v
\in \hunoc  \ \hbox{ a.e. in }\, (0,T)\,,
\end{aligned}
\\
& \label{inclu-z} {z(x,t) \in \ell(\teta_s(x,t)) \quad \forae\,
(x,t) \in \Gamma_c \times (0,T)\,,}
\\
 &
\begin{aligned}
 b({\bf {u}}_t,\vv)+a(\uu,\vv)+
\int_{\Omega} \teta \dive (\vv)& + \int_{\Gamma_c}(\chi{\bf u}
+\eeta )  \cdot{\bf v} \\ &= \pairing{\bfw'}{\bfw}{\mathbf{F}}{\vv}
 \quad \forall \vv\in \bfw \ \hbox{ a.e. in }\, (0,T)\,,
 \end{aligned}
\label{eqIa}\\
&\eeta\in \alpha(\uu) \ \ \text{in $(H^{-1/2}(\Gamma_c))^3 $}\ \
\hbox { a.e. in }\, (0,T), \label{incl1}
\\
 &\chi_t -\Delta\chi+ \xi
+\sig'(\chi)=-\lambda'(\chi) \teta_s-\frac 1 2\vert{\bf u}\vert^2
\quad\hbox{a.e. in } \Gamma_c\times
(0,T),\label{eqIIa}\\
 &\xi
\in \beta(\chi)\ \ \hbox { a.e. in }\, \Gamma_c\times (0,T),
\label{inclvincolo}\\
&\partial_{\nn_s}
  \chi=0\text{ a.e. in } \partial\Gamma_c\times (0,T)\,.\label{bordo1}
  \bigskip
\end{align}
\noindent Note that, to simplify notation, {we have incorporated}
the contribution $-\lambda'(\chi)\teta_{\rm eq}$ occurring
in~\eqref{eqII}  into the term $\sig'(\chi)$ in~\eqref{eqIIa}.
\begin{maintheorem}[Existence of a global solution]
 \label{t1}
Assume~\eqref{hyp:gamma-1}--\eqref{hypo-g}.
 Then, Problem~$(\mathbf{P})$
 admits at least a
global  solution $(  \teta, w, \teta_s,  z, \uu,\chi,\eeta,  \xi)$
with the regularity~\eqref{reg-teta}--\eqref{xireg}.
\end{maintheorem}
\begin{remark}[Positivity of the temperature]
\upshape {Clearly, in the case
\[
\dom(j) \subseteq (0,+\infty),
\]
(such as for~\eqref{scelta}),
 from relations~\eqref{inclu-w} and \eqref{inclu-z} we  infer that both
system temperatures $\teta$ and $\teta_s$ are strictly positive
almost everywhere in $\Omega \times (0,T)$ and in $\Gamma_c \times
(0,T)$, respectively.}
\end{remark}
\paragraph{Strategy of the proof of Theorem~\ref{t1}.}\,
We shall approximate Problem $(\mathbf{P})$ by suitably regularizing
equations \eqref{teta-weak} and \eqref{teta-s-weak}. {More
precisely, we shall add some viscosity terms to both equations, and
replace the operator $\ell$ therein by a Yosida-type regularization
(see \eqref{e:elle-astratto} below). For technical reasons (cf. with
Remark~\ref{rem:explanation}), we shall keep the viscosity parameter
distinct from the Yosida regularization parameter, denoting by
$\eps>0$ the former and by $\mu>0$ the latter. Hence, we shall call
$\pepsmu$ the initial and boundary value problem for the resulting
approximate system and prove that it is well-posed following this
outline: first} in Sections~\ref{s:3.1}--\ref{s:3.2} we are going to
prove the existence of a {\em local solution} by a fixed point
argument. Next, in Section~\ref{s:3.4} we are going to extend such a
solution to the whole interval $(0,T)$, while in Section~\ref{s:3.5}
we shall obtain contraction estimates leading to uniqueness for
Problem $\pepsmu$.

{Finally, in order to {prove}  Theorem~\ref{t1} we shall pass to the
limit in Problem $\pepsmu$ in two steps.  First, we shall keep
$\mu>0$ fixed and let $\eps \searrow 0$: in Section~\ref{s:4.1} we
are going to show that the {approximate solutions} converge, as
$\eps \searrow 0$, to a solution of the initial  and boundary value
problem obtained by setting $\eps=0$ in the approximate equations
\eqref{teta-weak-app}--\eqref{teta-s-weak-app} below. Secondly, we
shall also let $\mu \searrow 0$ and obtain in the limit a solution
to Problem $(\mathbf{P})$.}

\begin{notation}
\label{n:2} \upshape Henceforth,  for the sake of notational
simplicity, we shall use the same symbol $\langle \cdot, \cdot
\rangle$  for the duality pairings
 $\pairing{\bfw'}{\bfw}{\cdot}{\cdot}$,
$\pairing{V'}{V}{\cdot}{\cdot}$, and
$\pairing{{\hunoc}'}{\hunoc}{\cdot}{\cdot}$,  and, further,
 denote
by \begin{equation} \label{e:convention} \text{the symbols $C$, $C'$
most of the (positive) constants}
\end{equation}
occurring in calculations and estimates.
\end{notation}

\section{Approximation}
\label{s:3}
\subsection{The approximate problem}
\label{s:3.0}
  We  approximate Problem~$(\mathbf{P})$ by modifying
equations~\eqref{teta-weak} and~\eqref{teta-s-weak} in the following
way: {
\begin{itemize}
\item
first, we  shall add to~\eqref{teta-weak} the  regularizing
viscosity term $\eps \mathcal{R} (\teta_t)$  and
to~\eqref{teta-s-weak} the viscosity term $\eps
\mathcal{R}_{\Gamma_c} (\partial_t \teta_s)$, with $\eps>0$: this
shall  enable us to perform enhanced regularity estimates on the
(approximate) equations for the temperatures $\teta$
 and $\teta_s$ and ultimately to prove the global well-posedness of
 the approximate system, see  also Remark~\ref{rem:explanation};
\item
 second, both in~\eqref{teta-weak} and in~\eqref{teta-s-weak}  we shall replace the operator
$\ellv$  with its Yosida-type regularization
\begin{equation}
\label{e:elle-astratto} \calelleveps:= \left( \mu \Id +
\gammaveps\right)^{-1}\,,
\end{equation}
where $\mu>0$, $\Id: \R \to \R$ is the identity function, and
$\gammaveps: \R \to \R$ is the ($\mu$-)Yosida regularization of the
inverse $\gamma$ of $\ell$, see~\eqref{e:max-mon-2-v} below. The
choice of $\calelleveps$ is motivated by technical reasons, cf. with
Remarks~\ref{rem:explanation} and~\ref{rem:data} later on.
\end{itemize}
} \noindent
 {For the purposes of the above approximation, we recall
that the Yosida regularization  $\gammaveps: \R \to \R$ is a
Lipschitz continuous function, given by
\begin{equation}
\label{e:max-mon-2-v} \gammaveps(w):= \frac1\mu \left( w-\rhoveps
(w) \right)\qquad \text{for all $w\in \R$,}
\end{equation}
where $\rhoveps:\R \to \R$ is the ($\mu$-)resolvent operator
associated with $\gamma$, defined for every $w\in \R$
 as the unique solution
$\rhoveps (w) $ of the inclusion
\begin{equation}
\label{e:max-mon-1v}  \rhoveps (w)  -w + \mu \gamma(\rhoveps (w))
\ni 0.
\end{equation}
We shall also deal with the Yosida approximation $\jvstareps$ of
$\jvstar$, defined for every $\mu
>0$ by
\begin{equation}
\label{e:max-mon-4-bis}
\begin{gathered}
\jvstareps(w):= \min_{y \in \R} \left\{ \frac{|y-w|^2}{2\mu} +
\jvstar(w) \right\} \qquad \text{for all $w\in \R$.}
\end{gathered}
\end{equation}
We  recall that $\jvstareps \in C^1 (\R)$, with derivative
${\jvstareps}'=\gammaveps$, and that it
 fulfils
\begin{equation}
\label{e:max-mon-5-v}   \jvstareps(w) = \frac{\mu}2 | \gammaveps (w)
|^2 + \jvstar \left( \rhoveps(w) \right) \qquad \forall\, w \in
\R\,.
\end{equation}
 It was proved in~\cite[Sec.~3]{bfl1} that the operator
$\calelleveps$ is well-defined on $\R$, monotone, and that
\begin{equation}
\label{e:lip}\calelleveps  \colon \R \to \R \ \  \text{is  Lipschitz
continuous with Lipschitz constant $1/\mu$.}
\end{equation}
} {
\paragraph{Approximate initial data.}
In order to properly state our approximate problem, depending on the
parameters $\eps>0$ and $\mu>0$, we shall need to prescribe some
initial conditions for $\teta$ and $\teta_s$. To this aim, in
Lemma~\ref{dati-iniziali} we shall construct sequences of
 initial data $\tetazeroemu$ and $\tetaessezeroemu$ satisfying
\begin{align}
\label{h:teta-zero} &
 \tetazeroemu \in V\,, \qquad j^*_\mu \left( \calelleveps
 (\tetazeroemu)\right) \in L^1 (\Omega) \,,
\\
& \label{h:teta-s-zero}  \tetaessezeroemu \in \hunoc\,, \qquad
 j^*_\mu \left( \calelleveps
 (\tetaessezeroemu)\right) \in L^1 (\Gamma_c)\,,
\end{align}
and such that there exists a constant
 $M_0
>0$ \emph{independent of $\eps>0$ and $\mu>0$} with
\begin{align}
 \label{indip-eps-1}
 \left\|j^*_\mu \left( \calelleveps
 (\tetazeroemu)\right) \right\|_{L^1 (\Omega)} +  \left\|j^*_\mu \left( \calelleveps
 (\tetaessezeroemu)\right) \right\|_{L^1 (\Gamma_c)} \leq M_0 \quad
 \text{for all $\eps,\, \mu>0,$}
 \end{align}
and a constant $M_0^\mu
>0$ \emph{independent of $\eps>0$} but possibly depending on $\mu>0$
with
\begin{align}
 \label{indip-eps-2}
 \eps^{1/2} \|\tetazeroemu \|_{V}+
 \eps^{1/2} \|\tetaessezeroemu \|_{H^1 (\Gamma_c)}
\leq M_0^\mu \quad \text{for all $\eps>0$.}
\end{align}
For  simplicity, throughout this section we shall omit the indexes
$\eps$ and $\mu$ in the notation for the initial data.} The
variational formulation of the initial-boundary value problem
approximating Problem~($\mathbf{P}$) then reads:
\medskip\par\noindent
{\bf Problem $\pepsmu$}.
 Given a quadruple of  initial data
$(\tetazero,\tetaessezero, \uu_0, \chi_0)$
fulfilling~\eqref{h:teta-zero}--\eqref{indip-eps-2} and
\eqref{cond-uu-zero}--\eqref{cond-chi-zero},
 find  functions  $( \teta, \teta_s, \uu,\chi,
\xi, \eeta)$
 satisfying
\eqref{reguI}--\eqref{xireg},
 with $\teta$
and $\teta_s$ such that
\begin{align}
& \teta \in H^1 (0,T; V) \,, \label{reg-teta-forte}
\\
 & \teta_s \in  H^1 (0,T; \hunoc)\,, \label{reg-teta-s-forte}
\end{align}
fulfilling
\begin{align}
& \label{teta-weak-app}
\begin{aligned}
\eps \int_{\Omega} \teta_t \, v &+ \int_{\Omega}
\partial_t\calelleveps(\teta) v -\int_{\Omega} \dive(\uu_t) \, v
+\eps\int_{\Omega} \nabla \teta_t \, \nabla v +\int_{\Omega} \nabla
\teta \, \nabla v
\\ &+ \int_{\Gamma_c} k(\chi)
(\teta-\teta_s) v = \pairing{}{}{h}{v} \quad \forall\, v \in V \
\hbox{ a.e. in }\, (0,T)\,,
\end{aligned}
\\ & \label{teta-s-weak-app}
\begin{aligned}
&\eps \int_{\Gamma_c} \partial_t \teta_s \, v + \int_{\Gamma_c}
\partial_t\calelleveps(\teta_s) v
+\eps\int_{\Gamma_c}\nabla\partial_t{\teta_s}\nabla v
\\ &-\int_{\Gamma_c}
\partial_t \lambda(\chi) \, v    +\int_{\Gamma_c} \nabla \teta_s  \, \nabla
v = \int_{\Gamma_c} k(\chi) (\teta-\teta_s) v  \quad \forall\, v \in
\hunoc  \ \hbox{ a.e. in }\, (0,T)\,,
\end{aligned}
\end{align}
and equations~\eqref{eqIa}--\eqref{bordo1}, along with the initial
conditions \eqref{iniu} and \eqref{inichi} for $\uu$ and $\chi$,
respectively, and, for $\teta$ and $\teta_s$,
\begin{align}
& \teta(0)=\tetazero \qquad  \text{in $V$}\,, \label{teta-zero}
 \\
& \teta_s(0)=\tetaessezero \qquad  \text{in $\hunoc$}\,.
\label{teta-s-zero}
\end{align}
\begin{remark}
\label{e:furtherreg} \upshape Combining
\eqref{reg-teta-forte}--\eqref{reg-teta-s-forte} with the fact that
$\calelleveps$ is Lipschitz continuous (cf. \eqref{e:lip}), we
conclude that for any solution $( \teta, \teta_s, \uu,\chi, \xi,
\eeta)$ there holds
\begin{equation}
\label{e:wz-reg} \calelleveps(\teta) \in L^\infty (0,T;V) \cap H^1
(0,T;H), \qquad \calelleveps(\teta_s) \in L^\infty (0,T;H^1
(\Gamma_c)) \cap H^1 (0,T;L^2 (\Gamma_c))\,.
\end{equation}
\end{remark}
\begin{remark}
\label{rem:explanation} \upshape Without going into details, a few
comments on the approximate Problem $\pepsmu$ are in order. Besides
regularizing the maximal monotone operators in~\eqref{teta-weak} and
in~\eqref{teta-s-weak}, we have inserted in both equations a
viscosity term  in order to make each of (the Cauchy problems for)
the  approximate equations well-posed. {To understand why, let us
focus on equation~\eqref{teta-weak-app} (analogous considerations
apply to equation~\eqref{teta-s-weak-app}).}

 Indeed,
because of the nonlinear term $\partial_t (\calelleveps (\teta))$
in~\eqref{teta-weak-app}, in order to prove that the (Cauchy problem
for the) latter equation has a unique solution, one has to integrate
it in time. Hence, as it will be clear from the calculations in
Lemma~\ref{l:2}, the additional viscosity term $\eps
\mathcal{R}(\teta_t)$ enables us to deal, in the integrated version
of~\eqref{teta-weak-app},  with the third type boundary condition on
$\teta$. On the other hand, the choice of the operator
$\calelleveps$ in~\eqref{teta-weak-app} (of $\calelleps$
in~\eqref{teta-s-weak-app}), in place of the usual Yosida
regularization of $\ellv$, is due to technical reasons connected to
the construction of sequences of approximate initial data
fulfilling~\eqref{h:teta-zero}--\eqref{indip-eps-2}, see
Lemma~\ref{dati-iniziali}  {and Remark~\ref{rem:data}} later on.

 {The ultimate reason why we keep the viscosity parameter $\eps$ distinct from the Yosida parameter
 $\mu$ in both approximate equations \eqref{teta-weak} and \eqref{teta-s-weak}
is due to the fact that, in order to recover the
$L^\infty(0,T;H)$-regularity~\eqref{reg-w} for the solution
component $w$ (the $L^\infty
(0,T;L^2(\Gamma_c))$-regularity~\eqref{reg-zeta} for the solution
component $z$, respectively), one has to test some approximation of
\eqref{teta-weak} (\eqref{teta-s-weak}, respectively) by a function
approximating $w$ ($z$, resp.), and
 obtain some bound obviously independent of the approximation
parameter. In the present framework, performing such an estimate on
equation \eqref{teta-weak-app} (on \eqref{teta-s-weak-app}, resp.)
with $\eps>0$ would not lead to estimates on $\calelleveps (\teta)$
independent of the parameters $\eps$ and $\mu$, essentially because
the term $\pairing{}{}{\eps \mathcal{R}(\teta_t)}{\calelleveps
(\teta)}$ ($\pairing{}{}{\eps
\mathcal{R}_{\Gamma_c}(\partial_t\teta_s)}{\calelleveps (\teta_s)}$,
resp.) cannot be dealt with by monotonicity arguments. That is why,
in Section~\ref{s:4.2} we shall perform  the crucial estimate
leading to regularity \eqref{reg-w} and \eqref{reg-zeta} only after
taking the limit in Problem~$\pepsmu$ as  $\eps\searrow 0$ with $\mu
>0$ fixed.
 }
\end{remark}
\noindent The following theorem holds.
\begin{theorem}[Global well-posedness for Problem $\pepsmu$]
\label{t:3.1} Under assumptions~\eqref{hyp:gamma-1}--\eqref{hypo-g},
for any set of initial data $(\tetazero,  \tetaessezero, \uu_{0},
\chi_{0})$
 complying with
conditions~\eqref{h:teta-zero}--\eqref{indip-eps-2}
and~\eqref{cond-uu-zero}--\eqref{cond-chi-zero} and
 for all   $\eps, \ \mu
>0$
there exists a unique global solution
 $(  \teta, \teta_s, \uu,\chi, \xi, \eeta)$ to Problem~$\pepsmu$.
\end{theorem}
\noindent As already mentioned, we shall first of all prove the
existence of a local solution to Problem $\pepsmu$ by means of a
Schauder fixed point argument,  which relies on auxiliary
intermediate results on the existence and uniqueness of solutions
for the single equations of $\pepsmu$. We shall prove such results
in the following Section~\ref{s:3.1}, and conclude the proof of
local existence in Section~\ref{s:3.2}. Finally, in
Section~\ref{s:3.4} we shall show that the local solution extends to
a (unique, by Section~\ref{s:3.5})  global one.

\subsection{Fixed point setup}
\label{s:3.1} For a fixed $\mathsf{t} \in (0,T]$ (which shall be
specified later on) and a fixed constant $R>0$, we consider the~set
\begin{equation}\label{St}
\begin{aligned}
\mathcal{S}_{\mathsf{t}}:= \Big\{ (\teta,\teta_s, \chi) \in  & \,
L^2 (0,\mathsf{t}; H) \times L^{5/2} (0,\mathsf{t}; L^3 (\Gamma_c))
\times L^{10} (0,\mathsf{t}; L^6 (\Gamma_c)) \ : \\ &
 \| (\teta,\teta_s,
\chi)\|_{L^2 (0,\mathsf{t}; H) \times L^{5/2} (0,\mathsf{t}; L^3
(\Gamma_c)) \times L^{10} (0,\mathsf{t}; L^6 (\Gamma_c))} \leq R
\Big\}\,.
\end{aligned}
\end{equation}
 We shall construct an operator $\mathcal{T}$ {(its definition is split in
several steps)}, mapping $\mathcal{S}_{\widehat{T}}$ into itself for
a suitable time $0 < \widehat{T} \leq T$,  in such a way that any
fixed point of $\mathcal{T}$ yields a solution to Problem $\pepsmu$
on the interval $(0,\widehat{T})$. Then, in
Proposition~\ref{prop:loc} we shall prove that $\mathcal{T}:
\mathcal{S}_{\widehat{T}} \to \mathcal{S}_{\widehat{T}} $  admits a
fixed point by the Schauder theorem.
\begin{notation}
\label{n:3} \upshape We shall denote by
 {
$$
\begin{gathered}
\begin{gathered}
M_i, \ i=1, 2, 3 \ \ \text{a positive constant depending on $R$, on
the problem data,} \\
\text{ on $M_0$ (cf.~\eqref{indip-eps-1}) and $M_0^\mu$
(cf.~\eqref{indip-eps-2}),}
\\
\text{but  {\em independent} of $\eps>0$, and of the fixed $\,
\mathsf{t} \in (0,T]$}
\end{gathered}
\\
\begin{gathered}
M_i^\eps,\ i=1, 2, 3 \ \ \text{a positive constant depending on
$R$, on the problem data,} \\
\text{on $M_0$ (cf.~\eqref{indip-eps-1}) and $M_0^\mu$
(cf.~\eqref{indip-eps-2}), and possibly on $\eps>0$,}
\\
\text{
 but in any case {\em independent}  of  the fixed $\, \mathsf{t} \in
(0,T]$}\,. \end{gathered}
\end{gathered}
$$
} Further,  we shall keep to the notation~\eqref{e:convention} for
all the constants which do not depend on the approximating
parameter~$\varepsilon$.
\end{notation}
\paragraph{Step $1$.}
We take $(\widehat{\teta}, \widehat{\teta}_s, \widehat{\chi} ) \in
\mathcal{S}_{\mathsf{t}}$ and  consider (the Cauchy problem) for
system~\eqref{eqIa}--\eqref{bordo1}, in which $\teta$
in~\eqref{eqIa} is replaced by $\widehat{\teta}$ and $\chi$ and
$\teta_s$ on the right-hand side of~\eqref{eqIIa} are replaced by  $
\widehat{\chi} $ and $ \widehat{\teta}_s $, respectively.
\begin{lemma}
\label{l:1} Assume~\eqref{hyp:alpha}, \eqref{A5}, \eqref{hyp-sig},
\eqref{hyp-lambda}, \eqref{hypo-f},  \eqref{hypo-g},
\eqref{cond-uu-zero}, and \eqref{cond-chi-zero}.

 Then, there exists
a constant $M_1>0$ such that for all
 $(\widehat{\teta}, \widehat{\teta}_s, \widehat{\chi} ) \in
\mathcal{S}_{\mathsf{t}}$ there exists a unique quadruple $(\uu,
\eeta, \chi, \xi)$, with the regularity
\begin{align}
\label{stima-1} \| \uu \|_{H^1(0,\mathsf{t};\bfw)} &  +
\|\eeta\|_{L^2(0,\mathsf{t}; H^{-1/2}(\Gamma_c))} +
\|\xi\|_{L^2(0,\mathsf{t};L^2(\Gamma_c))}\\\nonumber &+\|\chi\|_{
L^{2}(0,\mathsf{t};H^2 (\Gamma_c))  \cap L^{\infty}(0,\mathsf{t};H^1
(\Gamma_c)) \cap H^1(0,\mathsf{t};L^2 (\Gamma_c))}\leq M_1\,,
\end{align}
complying with the initial conditions~\eqref{inichi}--\eqref{iniu},
solving the  PDE system
\begin{align}
\label{PDE1} &b({\bf {u}}_t,\vv)+a(\uu,\vv) +
\int_{\Gamma_c}(\chi{\bf u} +\eeta )  \cdot{\bf v} \\\nonumber &
=\pairing{}{}{\mathbf{F}}{\vv} - \int_{\Omega} \widehat{\teta} \dive
(\vv)
 \quad \forall \vv\in \bfw \ \hbox{ a.e. in }\, (0,\mathsf{t})\,,
  \\\nonumber
&\chi_t  -\Delta\chi+ \xi +\sig'(\chi) =-\lambda'(\widehat{\chi})
\widehat{\teta}_s-\frac 1 2\vert{\bf u}\vert^2 \quad\hbox{a.e. in }
\Gamma_c\times (0,\mathsf{t})\,,
\end{align}
and such that $\uu$ and $\eeta$ fulfil~\eqref{incl1}, $\chi$ and
$\xi$ comply with~\eqref{inclvincolo}--\eqref{bordo1} on
$(0,\mathsf{t})$.
\end{lemma}
\noindent
Let us point out that,  in view of~\eqref{hyp-lambda},  \eqref{St},
and the H\"{o}lder inequality, there holds
\begin{equation}
\label{stima-elle2} \| \lambda'(\widehat{\chi})\,
\widehat{\teta}_s\|_{L^2 (0,\mathsf{t}; L^2 (\Gamma_c))} \leq C \|
\widehat{\chi}\|_{L^{10} (0,\mathsf{t}; L^6 (\Gamma_c))}\, \|
\widehat{\teta}_s\|_{L^{5/2} (0,\mathsf{t}; L^3 (\Gamma_c))} \leq
CR^2\,.
\end{equation}
 Then, also taking into account the fact that
$\widehat \teta \in L^2 (0,\mathsf{t}; L^2 (\Gamma_c))$,
Lemma~\ref{l:1} follows from~\cite[Thm.~1]{bbr2}, to which we refer
the reader.

Now, we let
\begin{equation}\label{Vt}
\begin{aligned}
\mathcal{V}_{\mathsf{t}}:= \Big\{ (\uu,\chi) \in &\,
H^1(0,\mathsf{t};\bfw) \times  \big(L^{2}(0,\mathsf{t};H^2
(\Gamma_c)) \cap L^{\infty}(0,\mathsf{t};H^1 (\Gamma_c))  \cap
H^{1}(0,\mathsf{t};L^2 (\Gamma_c))\big) \ : \
\\ &
\|\uu\|_{H^1(0,\mathsf{t};\bfw) } + \|
\chi\|_{L^{2}(0,\mathsf{t};H^2 (\Gamma_c)) \cap
L^{\infty}(0,\mathsf{t};H^1 (\Gamma_c))  \cap H^{1}(0,\mathsf{t};L^2
(\Gamma_c)) } \leq M_1 \Big\}\,.
\end{aligned}
\end{equation}
Thanks to~Lemma~\ref{l:1}, we may define an operator
$$
\mathcal{T}_1 : \mathcal{S}_{\mathsf{t}} \to
\mathcal{V}_{\mathsf{t}}
$$
mapping every triple  $(\widehat{\teta}, \widehat{\teta}_s,
\widehat{\chi}) \in  \mathcal{S}_{\mathsf{t}}$ into the pair $(\uu,
\chi)$ solving the Cauchy problem for system~\eqref{PDE1} (we imply
that with $(\uu, \chi)$ the solution  components $\eeta$ and $\xi$
satisfying \eqref{incl1} and~\eqref{inclvincolo} are uniquely
associated).
\paragraph{Step $2$.}
We fix $(\widehat{\teta}_s, \overline{\uu}, \overline{\chi} ) \in
\pi_2 (\mathcal{S}_{\mathsf{t}}) \times \mathcal{V}_{\mathsf{t}}$
and consider (the Cauchy problem) for~\eqref{teta-weak-app} with
data~$(\widehat{\teta}_s, \overline{\uu}, \overline{\chi})$.
\begin{lemma}
\label{l:2} Assume~\eqref{hyp:gamma-1}, \eqref{hyp:gamma-2-a},
\eqref{hyp-k}, \eqref{hypo-h},
and~\eqref{h:teta-zero}--\eqref{indip-eps-2}.

Then, there exist $M_2, M_2^\eps>0$ such that  for all
 $(\widehat{\teta}_s, \overline{\uu}, \overline{\chi} ) \in \pi_2
(\mathcal{S}_{\mathsf{t}}) \times \mathcal{V}_{\mathsf{t}}$ there
exists a unique $\teta$ with
\begin{equation}
\label{stima-2} \| \teta \|_{L^2(0,\mathsf{t};V)\cap L^\infty
(0,\mathsf{t};L^1 (\Omega))}\leq M_2,
\quad\|\teta_t\|_{L^2(0,\mathsf{t};V)}\leq M_2^\eps\,,
\end{equation}
 complying with the
initial condition~\eqref{teta-zero} and with
\begin{equation}
\label{PDE2}
\begin{aligned}
\eps \int_{\Omega} \teta_t \, v
 +
\int_{\Omega}\partial_t\calelleveps(\teta) v
 & -\int_{\Omega} \dive(\overline{\uu}_t) \,
v  +\eps\int_{\Omega} \nabla \teta_t \, \nabla v +\int_{\Omega}
\nabla \teta \, \nabla v \\ &+ \int_{\Gamma_c} k(\overline{\chi})
(\teta-\widehat{\teta}_s) v
 = \pairing{}{}{h}{v}  \quad \forall\, v \in V \ \hbox{ a.e.
in }\, (0,\mathsf{t})\,.
\end{aligned}
\end{equation}
\end{lemma}
\noindent {\em Proof.} \ For simplicity, throughout the proof we
shall use the notation
\begin{equation}
\label{cal-K} \mathcal{K}(x,t):=k(\overline{\chi}(x,t)) \qquad
\forae\, (x,t) \in \Gamma_c \times (0,\mathsf{t})\,.
\end{equation}
 It follows from~\eqref{hyp-k},   the regularity
of $\overline{\chi}$ (cf. \eqref{Vt}),
and~\eqref{continuous-embedding-2} that
\begin{equation}
\label{e:4.2.1} \mathcal{K} \in L^\infty (0,\mathsf{t};H^1
(\Gamma_c))\,, \ \ \text{hence} \ \ \mathcal{K} \in L^\infty
(0,\mathsf{t};L^p (\Gamma_c)) \quad \forall\, 1 \leq p <\infty\,.
\end{equation}

In view of~\cite[Thm.~1]{dibenetto-showalter}, the Cauchy problem
for~\eqref{PDE2} has at least a solution $\teta \in
H^1(0,\mathsf{t};V)$. In order to prove uniqueness, we let
 $\teta_1, \, \teta_2 \in H^1 (0, \mathsf{t}; V)
$ be two solutions of the Cauchy problem
(\ref{teta-zero},~\ref{PDE2}), and set $\widetilde{\teta}:\teta_1
-\teta_2$. We subtract the equation for $\teta_2$ from the equation
for $\teta_1$ and integrate on  $(0,t)$, with $0 \leq t \leq
\mathsf{t}$. Thus, we get
$$
\begin{aligned}
\eps \int_{\Omega} \widetilde{\teta}(t) \, v + \int_{\Omega}
\left(\calelleveps(\teta_1 (t))- \calelleveps(\teta_2(t)) \right) v
 &+ \eps
\int_{\Omega} \nabla\widetilde{\teta}(t) \, \nabla v \\ & +
\int_{\Omega} (1 * \nabla \widetilde{\teta})(t) \, \nabla v  +
\int_{\Gamma_c} (1
* \mathcal{K} \widetilde{\teta})(t)\,  v =0
\end{aligned}
$$
for all $v \in V$. Hence, we take $v=\widetilde{\teta}$ and
integrate in time: also using that the operator $\calelleveps$ is
monotone, with straightforward computations we obtain
\begin{equation}
\label{e:4.2.2} \eps \int_0^t \|\widetilde{\teta}\|_H^2 +\eps
\int_0^t \|\nabla \widetilde{\teta}\|_H^2 +\frac12 \|(1 * \nabla
\widetilde{\teta})(t)\|_{H}^2  \leq \int_{0}^t \|\widetilde{\teta}
\|_{L^2(\Gamma_c)}\, \|(1 * \mathcal{K}\widetilde{\teta})
\|_{L^2(\Gamma_c)}.
\end{equation}
In view of the Young inequality for convolutions~\eqref{young}, we
have
$$
\|(1 * \mathcal{K}\widetilde{\teta})(s) \|_{L^2(\Gamma_c)} \leq
s^{1/2} \| \mathcal{K}\widetilde{\teta} \|_{L^2 (0,s;
L^2(\Gamma_c))}  \leq s^{1/2} \| \mathcal{K} \|_{L^\infty (0,s;
L^4(\Gamma_c))} \| \widetilde{\teta} \|_{L^2 (0,s; L^4(\Gamma_c))}
\quad \forall\, s \in [0,t]\,.
$$
Hence, \eqref{e:4.2.2} yields
\begin{equation}
\label{e:4.2.3}
 \eps \int_0^t \|\widetilde{\teta}\|_V^2 \leq \frac{\eps}2 \int_0^t \|\widetilde{\teta}\|_V^2
 + C_\eps t   \| \mathcal{K} \|_{L^\infty (0,\mathsf{t};
H^1(\Gamma_c))}^2 \int_0^t \|\widetilde{\teta} \|_{L^2 (0,s; V)}^2
\, \dd s\,,
\end{equation}
where $C_\eps$  depends on $\eps$ and also on the embedding
constants in
\eqref{continuous-embedding}--\eqref{continuous-embedding-2}. By
applying the Gronwall Lemma (see, e.g., \cite[Lemma A.3]{brezis73}),
we end up with
$$
\widetilde{\teta}(t) =0 \qquad \forae\, t \in (0, \mathsf{t})\,,
$$
whence the desired uniqueness.

We prove estimate~\eqref{stima-2} by testing~\eqref{PDE2} by
$\teta$.
 Using the
definition~\eqref{e:elle-astratto} of $\calelleveps$, we find
\begin{equation}
\label{e:chain-rule}
\begin{aligned}
\int_{\Omega}\partial_t\calelleveps(\teta) \teta &= \mu
\int_{\Omega}\partial_t\calelleveps(\teta)\calelleveps(\teta)
+\int_{\Omega}\partial_t\calelleveps(\teta)\gammaveps(\calelleveps(\teta))\\
& = \frac\mu2 \frac \dd{\dd t} \left\|\calelleveps(\teta)
\right\|_{H}^2 + \frac \dd{\dd t}\int_{\Omega}\jvstareps
\left(\calelleveps(\teta) \right),
\end{aligned}
\end{equation}
 the latter inequality ensuing from the chain rule for the convex
functional $\jvstareps$. Hence, upon integrating in time
 we get
 \begin{equation}
 \label{eqele1}
\begin{aligned}
\frac\varepsilon 2\|\teta(t)\|^2_V &+ \frac\mu2
\left\|\calelleveps(\teta(t)) \right\|_{H}^2 +
\int_{\Omega}\jvstareps \left(\calelleveps(\teta(t)) \right) +
\int_0^t \| \nabla\teta\|_{H}^2  + \int_0^t \int_{\Gamma_c}
\mathcal{K} \teta^2\\ & = \frac{\eps}2\|\tetazero\|^2_V+ \frac\mu2
\left\|\calelleveps(\tetazero) \right\|_{H}^2 +
\int_{\Omega}\jvstareps \left(\calelleveps(\tetazero) \right)+I_1 +
I_2+I_3\,,
\end{aligned}
\end{equation}
where the integral terms $I_i$, $i=1,2,3$ are specified
by~\eqref{est1}--\eqref{est3} below.  Now, in view
of~\eqref{hyp:gamma-2-a} (see Lemma~\ref{l:-ci-serve-anticipato}
later on), there exist $C_1,\,\overline{C}_2>0$ not depending on
$\varepsilon\in(0,1)$ such that
\begin{equation}
\label{coercivita} \mu \left\|\calelleveps(\teta(t)) \right\|_{H}^2
+\int_{\Omega} \jvstareps \left(\calelleveps(\teta(t)) \right) \geq
C_1 \| \teta(t) \|_{L^1 (\Omega)} -\overline{C}_2\,,
\end{equation}
whereas by~\eqref{h:teta-zero}
and~\eqref{indip-eps-1}--\eqref{indip-eps-2} we estimate the first
three summands on the right-hand side of~\eqref{eqele1}.
 Further, we estimate
\begin{equation}
\label{est1}
\begin{aligned}
 I_1 & =\int_0^t \int_{\Gamma_c} \mathcal{K}\,\widehat{\teta_s} \,
\teta \leq \int_0^t  \| \mathcal{K} \|_{L^6
(\Gamma_c)}\,\|\widehat{\teta_s} \|_{L^3 (\Gamma_c)}\, \|\teta
\|_{L^2 (\Gamma_c)}\\ &  \leq C \int_0^t  \| \mathcal{K} \|_{H^1
(\Gamma_c)}\,\|\widehat{\teta_s} \|_{L^3 (\Gamma_c)}\,
\left(\|\teta-m(\teta) \|_{V} + \| m(\teta) \|_{V}\right)
\\ & \leq C \int_0^t \| \mathcal{K} \|_{H^1
(\Gamma_c)}\,\|\widehat{\teta_s} \|_{L^3 (\Gamma_c)}\, \| \nabla
\teta \|_{L^2 (\Omega)} + C \int_0^t \| \mathcal{K} \|_{H^1
(\Gamma_c)}\,\|\widehat{\teta_s} \|_{L^3 (\Gamma_c)}\,\|
\teta\|_{L^1 (\Omega)} \\ &
\begin{aligned}
 \leq
 \rho \int_0^t  \|\nabla \teta
\|_{L^2 (\Omega)}^2 &+ C_{\rho}    \| \mathcal{K} \|_{L^{10}
(0,\mathsf{t}; H^1 (\Gamma_c))}^2\,\|\widehat{\teta_s}
\|_{L^{5/2}(0,\mathsf{t}; L^3 (\Gamma_c))}^2\\ &  +
 C \int_0^t \| \mathcal{K} \|_{H^1
(\Gamma_c)}\,\|\widehat{\teta_s} \|_{L^3 (\Gamma_c)}\,\|
\teta\|_{L^1 (\Omega)}\,,
\end{aligned}
\end{aligned}
\end{equation}
where the second inequality follows
from~\eqref{continuous-embedding}--\eqref{continuous-embedding-2},
the third one from the Poincar\'e inequality, and, thanks to Young's
inequality~\eqref{young-p-q}, the last one holds for a suitable
$\rho
>0$ to be chosen later. In the same way, we estimate
\begin{align}
& \label{est2}
\begin{aligned}
 I_2 & =\int_0^t \int_{\Omega}\dive(\overline{\uu}_t)\,\teta \leq
\int_0^t\|\dive(\overline{\uu}_t)\|_{H}\,\|\teta\|_{H}\\
&  \leq
C\int_0^t\|\dive(\overline{\uu}_t)\|_{H}\,\left(\|\teta-m(\teta)
\|_{V} + \| m(\teta) \|_{V}\right)
\\ &  \leq \rho \int_0^t  \| \nabla \teta
\|_{L^2 (\Omega)}^2 +{ C_{\rho}}' \int_0^t \| \overline{\uu}_t
\|_{\bfw}^2 + C\int_0^t  \| \overline{\uu}_t \|_{\bfw} \| \teta
\|_{L^1 (\Omega)},  \end{aligned}
\\
& \label{est3}
\begin{aligned}
 I_3  =\int_0^t\pairing{}{}{h}{\teta} &\leq \int_0^t
\| h \|_{V'} \, \|\teta\|_{V}\\ &  \leq \rho \int_0^t\|\nabla \teta
\|_{L^2 (\Omega)}^2 + {C_{\rho}}''\int_0^t\| h \|_{V'}^2 +
C\int_0^t\| h \|_{V'}\| \teta\|_{L^1 (\Omega)}\,.
\end{aligned}
\end{align}
Now, collecting~\eqref{eqele1}--\eqref{coercivita}
and~\eqref{est1}--\eqref{est3} (in which we choose $\rho\leq 1/6$),
taking into account~\eqref{e:4.2.1}, estimate~\eqref{stima-1} for
$\overline{\uu}$ and~\eqref{St} for $\widehat{\teta_s}$, using that
fourth term on the left-hand side of~\eqref{eqele1} is nonnegative
thanks to~\eqref{hyp-k}, and applying the Gronwall Lemma, we infer
that there exists $M_2
>0$, independent of $\eps$ and $\mu$, such that
\begin{equation}\label{ele1}
\varepsilon^{1/2}\|\teta(t)\|_V+\|\teta\|_{L^2(0,t;V)}+\|\teta(t)\|_{L^1(\Omega)}\leq
M_2 \qquad \forall\, t \in (0,\mathsf{t}].
\end{equation}

Secondly, we  test  \eqref{PDE2} by $\teta_t$.  Being $\calelleveps$
monotone, we easily see that
\begin{equation}
\label{e:positivity-bis} \int_{\Omega}\partial_t \calelleveps
(\teta)\teta_t \geq 0 \qquad \aein \, (0,T)\,.
\end{equation}
We now integrate in time: taking into
account~\eqref{e:positivity-bis},   using H\"older's inequality and
the Sobolev
embeddings~\eqref{continuous-embedding}--\eqref{continuous-embedding-2},
we easily conclude
\begin{equation}
\label{ele2}
\begin{aligned}
\eps\|\teta_t\|^2_{L^2(0,t;V)}+\frac 12\|\nabla\teta(t)\|^2_H \leq
\frac 1 2\|\nabla\tetazero\|^2_H &+ C
\int_0^t\|\mathcal{K}\|_{L^4(\Gamma_c)}\|\teta\|_{L^4(\Gamma_c)}\|\teta_t\|_{L^2(\Gamma_c)}\\
& +C
\int_0^t\|\mathcal{K}\|_{L^6(\Gamma_c)}\|\widehat\teta_s\|_{L^3(\Gamma_c)}\|\teta_t\|_{L^2(\Gamma_c)}
\\
& +C \int_0^t \| \dive(\overline{\uu}_t)\|_{H}\, \| \teta_t \|_{H}
+\int_0^t\|{h}\|_{V'}\|\teta_t\|_V
\\ &
\leq \frac12\| \teta_0 \|_{V}^2 + \frac12 \eps\|\teta_t\|^2_{L^2(0,t;V)}\\
& +C_{\eps}'
\int_0^t\|\mathcal{K}\|^2_{H^1(\Gamma_c)}\left(\|\teta\|^2_{L^4(\Gamma_c)}+\|\widehat\teta_s\|^2_{L^3(\Gamma_c)}\right)
\\ & +C_{\eps}{''} \left( \int_0^t \| \overline{\uu}_t \|_{\bfw}^2
+ \int_0^t \|h \|_{V'}^2\right).
\end{aligned}
\end{equation}
Eventually, in view of~\eqref{stima-1}, \eqref{St}, \eqref{e:4.2.1},
and \eqref{ele1}, from~\eqref{ele2} and the Poincar\'e inequality we
deduce
\begin{equation}\label{ele3}
\|\teta_t\|_{L^2(0,t;V)}+\|\teta(t)\|_{V}\leq M_2^\eps \qquad
\forall\, t \in (0,\mathsf{t}],
\end{equation}
from which \eqref{stima-2} follows. \fin

Due to  Lemma \ref{l:2},  we are in the position  of defining the
solution operator associated with~\eqref{PDE2}
\begin{equation}
\label{Wt}
 \mathcal{T}_2 : \pi_2 (\mathcal{S}_{\mathsf{t}}) \times
\mathcal{V}_{\mathsf{t}} \to \mathcal{W}_{\mathsf{t}}:=\left\{ \teta
\in H^1(0,\mathsf{t};V)\ : \ \| \teta \|_{L^2(0,\mathsf{t};V)\cap
L^\infty (0,\mathsf{t};L^1 (\Omega))} \leq
M_2\right\}\,.
\end{equation}

\paragraph{Step $3$.}
We fix $(\overline{\teta}, \widehat{\teta}_s,  \overline{\chi} ) \in
\mathcal{W}_{\mathsf{t}} \times  \pi_2 (\mathcal{S}_{\mathsf{t}})
\times \pi_2 (\mathcal{V}_{\mathsf{t}})$ and consider (the Cauchy
problem) for~\eqref{teta-s-weak-app} with data $(\overline{\teta},
\widehat{\teta}_s,  \overline{\chi} )$.
\begin{lemma}
\label{l:3} Assume~\eqref{hyp:gamma-1}, \eqref{hyp:gamma-2-a},
\eqref{hyp-k}, \eqref{hyp-lambda},  and \eqref{h:teta-s-zero}.

Then, there exist  $M_3,\,M_3^\eps>0$ such that  for all
$(\overline{\teta}, \widehat{\teta}_s, \overline{\chi} ) \in
\mathcal{W}_{\mathsf{t}} \times  \pi_2 (\mathcal{S}_{\mathsf{t}})
\times \pi_2 (\mathcal{V}_{\mathsf{t}})$ there exists a unique
$\teta_s \in H^1 (0,\mathsf{t};  \hunoc)$, with
\begin{equation}
\label{stima-3} \| \teta_s \|_{L^2 (0,\mathsf{t};  \hunoc)\cap
L^\infty (0,\mathsf{t};  L^1(\Gamma_c))} \leq M_3,\quad\|
\partial_t\teta_s \|_{L^2 (0,\mathsf{t};  \hunoc)} \leq M_3^\eps\,,
\end{equation}
such that $\teta_s$  complies with the initial
condition~\eqref{teta-s-zero},
 and
\begin{align}
\label{PDE3} &\eps \int_{\Gamma_c} \partial_t \teta_s \, v +
\int_{\Gamma_c}\partial_t\calelleps(\teta_s)v -\int_{\Gamma_c}
\partial_t \lambda(\overline{\chi}) \, v    +\int_{\Gamma_c} \nabla \teta_s  \, \nabla
v\\& +\eps \int_{\Gamma_c} \nabla \partial_t \teta_s \, \nabla v =
\int_{\Gamma_c} k(\overline{\chi})
(\overline{\teta}-\widehat{\teta}_s) v \quad \forall\, v \in \hunoc
\ \hbox{ a.e. in }\, (0,\mathsf{t})\,.
\end{align}
\end{lemma}
 \noindent
{\em Proof.} \, Thanks to~\cite[Thms.~1,4]{dibenetto-showalter}
there exists a unique solution $\teta_s \in H^1 (0,\mathsf{t};
\hunoc)$ to the Cauchy problem for~\eqref{PDE3}. Hence, we shall
just prove~\eqref{stima-3}, referring to  notation~\eqref{cal-K} for
the term $k(\overline{\chi})$.
 We proceed as for \eqref{stima-2}: hence, we   test \eqref{PDE3} by
 $\teta_s$ and integrate  in time. Developing the very same
 calculations as throughout
 \eqref{e:chain-rule}--\eqref{coercivita},
 exploiting~\eqref{hyp:gamma-2-a} (via Lemma~\ref{l:-ci-serve-anticipato}),
 and recalling \eqref{h:teta-s-zero} and \eqref{indip-eps-1}--\eqref{indip-eps-2}, we
  find
\[
\begin{aligned}
\frac\varepsilon 2\|\teta_s(t)\|^2_{\hunoc} & +C_1 \| \teta_s
(t)\|_{L^1(\Gamma_c)}
 + \int_0^t \|
\nabla\teta_s\|_{L^2 (\Gamma_c)}^2 \\ & \leq C + I_4 + I_5 + I_6\,,
\end{aligned}
\]
where, also in view~\eqref{hyp-lambda},
\begin{equation}
\label{i4} I_4=\int_0^t \int_{\Gamma_c} |\overline{\chi}_t|
|\lambda'(\overline{\chi})||\teta_s| \leq C\int_0^t
\|\overline{\chi}_t\|_{L^2 (\Gamma_c)}\, \left( \|
\overline{\chi}\|_{L^4 (\Gamma_c)}+1 \right)\, \| \teta_s\|_{L^4
(\Gamma_c)}\,,
\end{equation}
and, in view of~\eqref{e:4.2.1},
\begin{align}
& \label{i5} I_5 \leq   \int_0^t  \|\mathcal{K}\|_{L^4 (\Gamma_c)}
\, \| \overline{\teta}\|_{L^4 (\Gamma_c)}\, \| \teta_s\|_{L^2
(\Gamma_c)} \leq \| \mathcal{K}\|_{L^\infty (0,t;L^4 (\Gamma_c))}
\int_0^t   \| \overline{\teta}\|_{L^4 (\Gamma_c)}\, \|
\teta_s\|_{L^2 (\Gamma_c)}\,,
\\ &
\label{i6}
\begin{aligned}
I_6 &=\int_0^t \int_{\Gamma_c} \|\mathcal{K}\|_{L^6 (\Gamma_c)} \,
\| \widehat{\teta}_s\|_{L^3 (\Gamma_c)}\, \|\teta_s \|_{L^2
(\Gamma_c)}
\\ & \leq \| \mathcal{K}\|_{L^\infty (0,t;L^6 (\Gamma_c))} \int_0^t
\| \widehat{\teta}_s\|_{L^3 (\Gamma_c)}\, \|\teta_s \|_{L^2
(\Gamma_c)}\,.
\end{aligned}
\end{align}
 Taking into account~\eqref{St},
 \eqref{Vt}, and
 \eqref{Wt}, and estimating  the term $\| \teta_s \|_{L^2
 (\Gamma_c)}$ by $\| \teta_s \|_{L^1
 (\Gamma_c)}$ and $\| \nabla \teta_s \|_{L^2 (\Gamma_c)}^2$
  in the  same
  as in~\eqref{est1}--\eqref{est2}, we finally apply the
 Gronwall Lemma to conclude
\begin{equation}\label{stimaM3}
\eps^{1/2}\|\teta_s\|_{L^\infty(0,\mathsf{t};\hunoc)}+\|\teta_s\|_{L^2(0,\mathsf{t};H^1(\Gamma_c))
\cap L^\infty(0,\mathsf{t};L^1(\Gamma_c))}\leq M_3.
\end{equation}

Then,  we test~\eqref{PDE3} by $\partial_t{\teta_s}$ and integrate
in time. Thanks to~\eqref{stimaM3} and arguing in the same way as
in~\eqref{i4}--\eqref{i5}, we find
$$
\begin{aligned}
\eps\int_0^t \|\partial_t{\teta_s}\|_{\hunoc}^2
&+\frac12\|\teta_s(t)\|_{H^1(\Gamma_c)}^2 +\int_0^t
\int_{\Gamma_c}\partial_t \calelleps (\teta_s)\partial_t \teta_s
\\ &
\begin{aligned}
\leq \frac12\|\tetaessezero\|_{H^1(\Gamma_c)}^2 &+ \int_0^t
\|\overline{\chi}_t\|_{L^2 (\Gamma_c)}\, \| \overline{\chi}\|_{L^4
(\Gamma_c)}\, \| \partial_t\teta_s\|_{L^4 (\Gamma_c)} \\ & +\int_0^t
\|\mathcal{K}\|_{L^4 (\Gamma_c)}  \, \| \overline{\teta}\|_{L^4
(\Gamma_c)}\, \|
\partial_t \teta_s\|_{L^2 (\Gamma_c)} \\ & + \int_0^t \|\mathcal{K}\|_{L^6 (\Gamma_c)}  \, \| \widehat{\teta}_s\|_{L^3
(\Gamma_c)}\, \|
\partial_t \teta_s\|_{L^2 (\Gamma_c)}\,.
\end{aligned}
\end{aligned}
$$
Again taking into account~\eqref{St}, \eqref{Vt}, and~\eqref{Wt}, as
well as the fact that the third summand on the left-hand side is
nonnegative, we readily deduce the second part of \eqref{stima-3}.
\fin
\\
  Lemma~\ref{l:3}
enables us to define a solution operator associated
with~\eqref{PDE3}
$$ \mathcal{T}_3 : \mathcal{W}_{\mathsf{t}} \times  \pi_2
(\mathcal{S}_{\mathsf{t}}) \times \pi_2 (\mathcal{V}_{\mathsf{t}})
\to \mathcal{Y}_{\mathsf{t}}:= \left\{ \teta_s \in  L^\infty
(0,\mathsf{t}; \hunoc) \ : \ \|\teta_s \|_{L^2 (0,\mathsf{t};
\hunoc)\cap L^\infty (0,\mathsf{t}; L^1(\Gamma_c))} \leq M_3
\right\}\,.
$$

%
%

%

\subsection{Local existence for Problem $\pepsmu$}
\label{s:3.2}

\begin{proposition}
\label{prop:loc}
 Assume \eqref{hyp:gamma-1}--\eqref{hypo-g},
\eqref{cond-uu-zero}--\eqref{cond-chi-zero}
and~\eqref{h:teta-zero}--\eqref{indip-eps-2}. Then,  there exists
$\widehat{T} \in (0,T]$, {possibly depending on $\mu>0$,} such that
for every $\eps>0$
 Problem $\pepsmu$ admits a solution $(  \teta,
\teta_s, \uu,\chi,\xi,\eeta)$ on the interval $(0,\widehat{T})$.
\end{proposition}

{\em Proof.} \ In view of the above Lemmata \ref{l:1}, \ref{l:2},
\ref{l:3}  we are able to define an operator $\mathcal{T}$ whose
fixed points are solutions of Problem $\pepsmu$.

\paragraph{Definition of $\mathcal{T}$.} In the end, we define
$$
\mathcal{T} : \mathcal{S}_{\mathsf{t}} \to \mathcal{W}_{\mathsf{t}}
\times \mathcal{Y}_{\mathsf{t}} \times \pi_2
(\mathcal{V}_{\mathsf{t}})
$$
by setting for every $(\widehat{\teta}, \widehat{\teta}_s,
\widehat{\chi} ) \in   \mathcal{S}_{\mathsf{t}}$
\begin{equation}
\label{def-T} \mathcal{T}(\widehat{\teta}, \widehat{\teta}_s,
\widehat{\chi} ) := (\teta, \teta_s,\chi), \ \ \text{where} \ \
\begin{cases}
\teta= \mathcal{T}_2 ( \widehat{\teta}_s\,, \mathcal{T}_1
(\widehat{\teta}, \widehat{\teta}_s, \widehat{\chi} ))\,, \\
\teta_s = \mathcal{T}_3 (\mathcal{T}_2 ( \widehat{\teta}_s,\,
\mathcal{T}_1 (\widehat{\teta}, \widehat{\teta}_s, \widehat{\chi}
)), \, \widehat{\teta}_s, \, \pi_2 (\mathcal{T}_1 (\widehat{\teta},
\widehat{\teta}_s, \widehat{\chi} )))\,,
\\
\chi=\pi_2 (\mathcal{T}_1 (\widehat{\teta}, \widehat{\teta}_s,
\widehat{\chi} ))\,.
\end{cases}
\end{equation}
Thus, in order to prove Proposition \ref{prop:loc} it is sufficient
to show that
 there exists $\widehat{T} \in (0,T]$ such that for every $\eps >0$
\begin{gather}
\label{itself} \text{
 $\mathcal{T}$ maps
 $\mathcal{S}_{\widehat{T}}$ into itself,}
 \\
 \label{compact-conti}
 \begin{gathered}
\mathcal{T} : \mathcal{S}_{\widehat{T}} \to
\mathcal{S}_{\widehat{T}}  \ \ \text{ is compact and
 continuous} \\ \text{ w.r.t. the topology of $L^2 (0,\widehat{T} ; H) \times L^{5/2} (0,\widehat{T} ; L^3 (\Gamma_c))
\times L^{10} (0,\widehat{T} ; L^6 (\Gamma_c))$.}
\end{gathered}
\end{gather}
\paragraph{Ad~\eqref{itself}.}\, We fix $(\widehat{\teta}, \widehat{\teta}_s,
\widehat{\chi} ) \in   \mathcal{S}_{\mathsf{t}}$ and let
$(\teta,\teta_s,\chi)=\mathcal{T}(\widehat{\teta},
\widehat{\teta}_s, \widehat{\chi} )$.
 First of all,
the three-dimensional version of the Gagliardo-Nirenberg inequality
(cf.~\cite[p.~125]{nier}) yields
\begin{equation}
 \|v\|_H\leq C\|\teta\|^{3/5}_V\|\teta\|^{2/5}_{L^1(\Omega)},
\end{equation}
so that,
 using \eqref{stima-2}, we  get
\begin{equation}
\label{dentro-la-palla-1} \|\teta\|_{L^2(0,\mathsf{t};H)}\leq \|
\teta\|_{L^2(0,\mathsf{t};V)}^{3/5}\, \|
\teta\|_{L^2(0,\mathsf{t};L^1(\Omega))}^{2/5} \leq
\mathsf{t}^{1/5}M_2\,.
\end{equation}
 The Gagliardo-Nirenberg
inequality in $2D$ gives
$$
\| \teta_s \|_{L^3(\Gamma_c)} \leq C
\|\teta_s\|^{2/3}_{\hunoc}\|\teta_s\|^{1/3}_{L^1(\Gamma_c)},
$$
whence
\begin{equation}
\label{dentro-la-palla-2}
\|\teta_s\|_{L^{5/2}(0,\mathsf{t};L^3(\Gamma_c))}\leq \|
\teta\|_{L^2(0,\mathsf{t};\hunoc)}^{2/3}\, \|
\teta\|_{L^5(0,\mathsf{t};L^1(\Gamma_c))}^{1/3} \leq
\mathsf{t}^{1/{15}}M_3\,.
\end{equation}
Finally, by~\eqref{stima-1} we have
\begin{equation}
 \|\chi\|_{L^{10}(0,\mathsf{t}); L^6(\Gamma_c)}\leq C
 \mathsf{t}^{1/{10}}M_1\,,
\end{equation}
the constant $C$ depending on the Sobolev
embedding~\eqref{continuous-embedding-2}.
 Clearly, there exists $\widehat{T}$ (which does not depend on $\eps$) such that
$(\teta,\teta_s,\chi)$ belongs to ${\cal S}_{\widehat{T}}$, hence
the operator ${\cal T}$ maps ${\cal S}_{\widehat T}$ into itself.
\paragraph{Ad~\eqref{compact-conti}.}
Exploiting \eqref{stima-1}, \eqref{stima-2}, and \eqref{stima-3},
the Sobolev
embeddings~\eqref{continuous-embedding}--\eqref{continuous-embedding-2}
and~\cite[Thm.~4, Cor.~5]{Simon87}, one sees immediately that the
operator $\mathcal{T}: \mathcal{S}_{\widehat{T}} \to
\mathcal{S}_{\widehat{T}}$ is compact. We shall prove that
$\mathcal{T}$ is continuous in three steps, basically  checking that
the operators $\mathcal{T}_i$, $i=1, 2, 3$ defined in
Section~\ref{s:3.1} are continuous w.r.t. to suitable topologies.

 We fix a sequence $\{(\widehat{\teta}_n,
\widehat{\teta}_{s,n}, \widehat{\chi}_n ) \} \subset
\mathcal{S}_{\widehat{T}}$ such that there exists
$(\widehat{\teta}_\infty,  \widehat{\teta}_{s,\infty},
\widehat{\chi}_\infty )  \in  \mathcal{S}_{\widehat{T}}$ with
\begin{align}
& \label{co:1} \widehat{\teta}_n \to \widehat{\teta}_\infty \qquad
\text{in $L^2 (0,\widehat{T}; H)$} \quad \text{as $n \to \infty$\,,}
\\
& \label{co:2} \widehat{\teta}_{s,n} \to \widehat{\teta}_{s,\infty}
\qquad \text{in $L^{5/2} (0,\widehat{T}; L^3 (\Gamma_c))$} \quad
\text{as $n \to \infty$\,,}
\\
& \label{co:3} \widehat{\chi}_{n} \to \widehat{\chi}_{\infty} \qquad
\text{in $L^{10} (0,\widehat{T}; L^6 (\Gamma_c))$} \quad \text{as $n
\to \infty$\,,}
\end{align}
we let $(\uu_n, \chi_n):= \mathcal{T}_1 (\widehat{\teta}_n,
\widehat{\teta}_{s,n}, \widehat{\chi}_n )$,  and denote by $\{
\eeta_n \}  $ and $\xi_n$ the associated  sequences of selections of
the graph $\alpha$ and $\beta$, respectively, such that
\eqref{incl1} and \eqref{inclvincolo} hold for all $n \in \N$. Due
to~\eqref{stima-1}, there exists a subsequence (which we do not
relabel) and a quadruple
$$
\begin{aligned}
(\uu_\infty, \eeta_\infty, \chi_\infty, \xi_\infty) \in
H^1(0,\widehat{T};\bfw) & \times L^2(0,\widehat{T};
H^{-1/2}(\Gamma_c))\\ & \times (L^{2}(0,\widehat{T};H^2 (\Gamma_c))
\cap L^\infty(0,\widehat{T};H^1 (\Gamma_c)) \cap H^1
(0,\widehat{T};L^2 (\Gamma_c)))\\ & \times L^2(0,\widehat
T;L^2(\Gamma_c))
\end{aligned}
$$
such that  the following convergences hold as $n \nearrow \infty$:
\begin{align}
& \label{co:4}
\begin{aligned}
& \uu_{n}  \weakto  \uu_{\infty} \ \  \text{ in } \
H^1(0,\widehat{T}; \bfw),
 \qquad
 \uu_{n} \to \uu_{\infty}  \ \ \text{ in } \ C^0([0,\widehat{T}];
(H^{1-s}(\Omega))^3)\ \ \text{for all $s>0$,}
\\
 &\uu_{n}  \to \uu_\infty \ \text{ in } \ C^0([0,\widehat{T}];
 (L^p(\Gamma_c))^3), \quad
 \text{for all $1 \leq p <4$,}
 \end{aligned}
 \\
 &
 \label{co:5}
\eeta_n \weakto \eeta_\infty \ \
 \text{ in } \
L^2(0,\widehat{T};(H^{-1/2}(\Gamma_c))^3)\,,
\\
& \label{co:6}
\begin{aligned}
& \chi_n   \weaksto \chi_\infty  \ \  \text{ in } \
H^1(0,\widehat{T};L^2(\Gamma_c)) \cap
L^\infty(0,\widehat{T};H^1(\Gamma_c))\cap
L^{2}(0,\widehat{T};H^2(\Gamma_c)),
\\
 &\chi_{n}  \to \chi_\infty  \  \ \text{ in } \
C^0([0,\widehat{T}];H^{1-s}(\Gamma_c))\cap L^2(0,\widehat
T;H^{2-s}(\Gamma_c)) \ \ \text{for all $s>0$}\,,
\end{aligned}
\\
& \label{co:6bis} \xi_n\weakto\xi_\infty \ \ \text{ in
}L^2(0,\widehat T;L^2(\Gamma_c))\,.
\end{align}
Thanks to  these convergences, and recalling \eqref{hyp-lambda}, it
is possible to pass to the limit
 in the
PDE system~\eqref{PDE1} as $n \nearrow \infty$ (see the proof
of~\cite[Prop.~4.7]{bbr1} for further details). Thus, we   prove
that the quadruple $(\uu_\infty, \eeta_\infty,
\chi_\infty,\xi_\infty)$ fulfils~\eqref{PDE1},  with data
$(\widehat{\teta}_\infty, \widehat{\teta}_{s,\infty},
\widehat{\chi}_\infty )$, on $(0,\widehat{T})$. Furthermore, the
pair $(\uu_\infty, \chi_\infty)$ clearly complies with  the initial
conditions~\eqref{iniu}--\eqref{inichi}. Note that the
identifications
$$
\xi_\infty\in\beta(\chi_\infty),\quad\eeta_\infty\in\alpha(\uu_\infty)
$$
are proved by semicontinuity arguments (see \cite[Lemma~1.3,
p~42]{barbu}). In the end,  recalling the definition of the operator
$\mathcal{T}_1 $ we conclude that
\begin{equation}
\label{ope-t1}
 (\uu_\infty, \chi_\infty)=\mathcal{T}_1
(\widehat{\teta}_\infty,  \widehat{\teta}_{s,\infty},
\widehat{\chi}_\infty )\,.
\end{equation}
In fact,  by uniqueness of the limit we readily deduce that
convergences~\eqref{co:4}--\eqref{co:6bis} hold for  the whole
sequences. In particular, we have  checked
that~\eqref{co:1}--\eqref{co:3} imply
\begin{equation}
\label{ded-1} \mathcal{T}_1 (\widehat{\teta}_n,
\widehat{\teta}_{s,n}, \widehat{\chi}_n ) \to \mathcal{T}_1
(\widehat{\teta}_\infty, \widehat{\teta}_{s,\infty},
\widehat{\chi}_\infty ) \qquad \text{in the sense of \eqref{co:4},
\eqref{co:6}. }
\end{equation}

 We now
consider the sequence $\teta_n: =\mathcal{T}_2
(\widehat{\teta}_{s,n}, \uu_n, \chi_n) =\mathcal{T}_2
(\widehat{\teta}_{s,n},
  \mathcal{T}_1 (\widehat{\teta}_n,
\widehat{\teta}_{s,n}, \widehat{\chi}_n ))$ for all $n \in \N$.
Thanks to~\eqref{stima-2}, $\{ \teta_n \}$ is bounded in $H^1
(0,\widehat{T}, V)$, hence there exists $\teta_\infty \in H^1
(0,\widehat{T}, V)$ such that (up to  a subsequence)
\begin{equation}
\label{co:8}
\begin{aligned}
&
 \teta_{n}  \weakto  \teta_{\infty} \ \  \text{ in } \ H^1(0,\widehat{T}; V),
 \qquad
 \teta_{n} \to \teta_{\infty}  \ \ \text{ in } \ C^0([0,\widehat{T}];
H^{1-s}(\Omega))\ \ \forall\, s>0,
\\
 &\teta_{n}  \to \teta_\infty \ \text{ in } \ C^0([0,\widehat{T}];
 L^p(\Gamma_c)) \quad
 \text{for all $1 \leq p <4$,}
\end{aligned}
\end{equation}
as $n \nearrow \infty$. In particular, $\teta_\infty$ complies
with~\eqref{teta-zero}. Moreover, using that $\calelleveps$ is
Lipschitz continuous,  {taking into account the strong convergence
for $\teta_n$  specified by the second of \eqref{co:8} and recalling
that, by maximal monotonicity, the graph of $\calelleveps$ is
strongly-weakly closed, it is not difficult to conclude that
\begin{equation}
\label{co:9}
\begin{cases}
&\calelleveps(\teta_n) \weaksto \calelleveps(\teta_\infty) \qquad
\text{in $L^\infty (0,\widehat{T};V) \cap H^1 (0,\widehat{T};H),$}
\\
& \calelleveps(\teta_n) \to \calelleveps(\teta_\infty) \qquad
\text{in $ C^0([0,\widehat{T}]; H)$.} \end{cases}
\end{equation}
} Furthermore, it follows from \eqref{co:2}, \eqref{co:6},
\eqref{co:8},  and the Lipschitz continuity of $k$ that, among
others, the following convergences hold $n \nearrow \infty$:
\begin{equation}
\label{co:10}
\begin{gathered}
k(\chi_n) \teta_n \to k(\chi_\infty) \teta_\infty \quad \text{in
$L^\infty(0,\widehat{T};  L^2(\Gamma_c)) $ for all $1 \leq p <4$, and} \\
k(\chi_n) \widehat{\teta}_{s,n} \to k(\chi_\infty)
\widehat{\teta}_{s,\infty} \quad \text{in $ L^2(0,\widehat{T};
L^2(\Gamma_c))$}\,.
\end{gathered}
\end{equation}
Combining  \eqref{co:4} with \eqref{co:8}--\eqref{co:10}, we   pass
to the limit as $n \nearrow \infty$ in~\eqref{PDE2} with data
$(\widehat{\teta}_{s,n}, \uu_n, \chi_n)$, and we deduce that
$\teta_\infty$ solves (the Cauchy problem) for equation
\eqref{PDE2},  with the triple
 $(\widehat{\teta}_{s,\infty}, \uu_\infty,
\chi_\infty)$, on $(0,\widehat{T})$. Thus,
$$
\teta_\infty= \mathcal{T}_2 (\widehat{\teta}_{s,\infty}, \uu_\infty,
\chi_\infty)= \mathcal{T}_2
(\widehat{\teta}_{s,\infty},\mathcal{T}_1 (\widehat{\teta}_\infty,
\widehat{\teta}_{s,\infty}, \widehat{\chi}_\infty ))\,,
$$
the second equality ensuing from \eqref{ope-t1}. Again, since the
limit $\teta_\infty$ does not depend on the  subsequence in
\eqref{co:8}, it turns out  the convergences specified therein hold
along the whole sequence $\{ \teta_n \}$. In conclusion,
\begin{equation}
\label{ded-2} \mathcal{T}_2 (\widehat{\teta}_{s,n}, \mathcal{T}_1
(\widehat{\teta}_n, \widehat{\teta}_{s,n}, \widehat{\chi}_n ) ) \to
 \mathcal{T}_2 (\widehat{\teta}_{s,\infty},\mathcal{T}_1
(\widehat{\teta}_\infty, \widehat{\teta}_{s,\infty},
\widehat{\chi}_\infty )) \qquad \text{in the sense of \eqref{co:8}.
}
\end{equation}

Finally, we let $$\teta_{s,n}:= \mathcal{T}_3 (\teta_n,
\widehat{\teta}_{s,n}, \chi_n) = \mathcal{T}_3 (\mathcal{T}_2 (
\widehat{\teta}_{s,n},\, \mathcal{T}_1 (\widehat{\teta}_n,
\widehat{\teta}_{s,n}, \widehat{\chi}_n )), \,
\widehat{\teta}_{s,n}, \, \pi_2 (\mathcal{T}_1 (\widehat{\teta}_n,
\widehat{\teta}_{s,n}, \widehat{\chi}_n )))$$ for every $n \in \N$.
By \eqref{stima-3}, $\{\teta_{s,n}\}$ is bounded in $H^1
(0,\widehat{T}; \hunoc)$. Thus,
  there exists a (not relabeled)
subsequence and $\teta_{s,\infty}\in H^1 (0,\widehat{T}; \hunoc )$
such that {
\begin{align}
& \label{co:11}
\begin{aligned}
&
 \teta_{s,n} \weakto \teta_{s,\infty} \qquad \text{in $H^1
(0,\widehat{T}; H^1(\Gamma_c) )$,}
\\
&
 \teta_{s,n} \to \teta_{s,\infty} \qquad \text{in $C^0
([0,\widehat{T}]; H^{1-s}(\Gamma_c) )$ for all $s >0$,}
\end{aligned}
\\
& \label{co:12}
\begin{aligned}
&
 \calelleps(\teta_{s,n}) \weaksto \calelleps(\teta_{s,\infty}) \qquad\text{in $L^\infty (0,\widehat{T};
 \hunoc) \cap
 H^1 (0,\widehat{T};  L^2(\Gamma_c) )$,}
 \\
&
 \calelleps(\teta_{s,n}) \to \calelleps(\teta_{s,\infty}) \qquad\text{in $C^0 ([0,\widehat{T}]; L^2(\Gamma_c) )$,}
\end{aligned}
\end{align}
where the latter convergences can be proved arguing in the very same
way as for \eqref{co:9}.} In particular, $\teta_{s,\infty} $ fulfils
initial condition \eqref{teta-s-zero}, and, like in the proof of
Lemma \ref{l:3}, combining convergences \eqref{co:11}--\eqref{co:12}
we conclude that $\teta_{s,\infty}$ fulfils \eqref{PDE3}  on
$(0,\widehat{T})$, with data
 $(\teta_\infty,\widehat{\teta}_{s,\infty},
\chi_\infty)$.
 Now,  \eqref{co:6} and the growth
properties of $\lambda$ (see \eqref{hyp-lambda}) clearly yield that
\begin{equation}
\label{co:13} \lambda' (\chi_n)\partial_t\chi_n\weak  \lambda'
(\chi_\infty)\partial_t\chi_\infty \quad \text{in $L^2
(0,\widehat{T}; \hunoc' )$.}
\end{equation}
 Then, exploiting \eqref{co:6} and \eqref{co:11} we easily check
that
\begin{equation}
\label{co:14} k(\chi_n) \teta_{s,n} \to k(\chi_\infty)
\teta_{s,\infty} \qquad \text{in $L^\infty (0,\widehat{T};
L^{2}(\Gamma_c) )$.}
\end{equation}
Collecting \eqref{co:11}--\eqref{co:14} and also taking into account
\eqref{co:10}, we pass to the limit in \eqref{PDE3} (with data
$(\teta_n, \widehat{\teta}_{s,n}, \chi_n)$) as $n \nearrow \infty$.
Hence,
$$
\teta_{s,\infty} =\mathcal{T}_3 (\teta_\infty,\widehat{\teta}_{s,
\infty}, \chi_\infty )= \mathcal{T}_3 (\mathcal{T}_2
(\widehat{\teta}_{s,\infty},\,\mathcal{T}_1 (\widehat{\teta}_\infty,
\widehat{\teta}_{s,\infty}, \widehat{\chi}_\infty )),\,
\widehat{\teta}_{s, \infty},\, \pi_2 (\mathcal{T}_1
(\widehat{\teta}_\infty, \widehat{\teta}_{s,\infty},
\widehat{\chi}_\infty ))),
$$
and we again deduce that convergences \eqref{co:11}--\eqref{co:12}
hold for the whole sequences  $\{ \teta_{s,n} \}$. Eventually, we
have that
\begin{equation}
\label{ded:3}
\begin{aligned}
& \mathcal{T}_3 (\mathcal{T}_2 (\widehat{\teta}_{s,n},\,
\mathcal{T}_1 (\widehat{\teta}_n, \widehat{\teta}_{s,n},
\widehat{\chi}_n )),\, \widehat{\teta}_{s, n},\, \pi_2
(\mathcal{T}_1 (\widehat{\teta}_n, \widehat{\teta}_{s,n},
\widehat{\chi}_n ))) \, \to \\ & \mathcal{T}_3 (\mathcal{T}_2
(\widehat{\teta}_{s,\infty},\,\mathcal{T}_1 (\widehat{\teta}_\infty,
\widehat{\teta}_{s,\infty}, \widehat{\chi}_\infty )),\,
\widehat{\teta}_{s, \infty},\, \pi_2 (\mathcal{T}_1
(\widehat{\teta}_\infty, \widehat{\teta}_{s,\infty},
\widehat{\chi}_\infty ))) \quad \text{in the sense of
\eqref{co:11}.}
\end{aligned}
\end{equation} Clearly, \eqref{ded-1}, \eqref{ded-2}
and \eqref{ded:3} show  that $\mathcal{T}$ is continuous in the
sense of \eqref{compact-conti}. \fin
\subsection{Global existence for Problem $\pepsmu$}
\label{s:3.4}
In order to  show that the local solution to  Problem $\pepsmu$
actually extends to the whole time interval $(0,T)$, we shall obtain
some global in time estimates on the solution components
$(\teta,\teta_s,\uu,\chi)$ and then   use a fairly standard argument
to conclude that, for every $\eps>0$ and $\mu>0$, the local solution
found in Proposition~\ref{prop:loc} extends to the (unique, by the
calculations in Section~\ref{s:3.5}) global solution of Problem
$\pepsmu$.

\begin{lemma}[Global estimates]
\label{l:global-esti} Assume~\eqref{hyp:gamma-1}--\eqref{hypo-g} and
let  $(\tetazero,  \tetaessezero, \uu_0, \chi_0)$ be a quadruple of
initial data
 complying with
conditions~\eqref{h:teta-zero}--\eqref{indip-eps-2}
and~\eqref{cond-uu-zero}--\eqref{cond-chi-zero}. Then, for every
$\mu>0$ there exists a constant $C_{\mu} >0$, depending on the
problem data, on $M_0$ (cf.~\eqref{indip-eps-1}) {and $M_0^{\mu}$
(cf.~\eqref{indip-eps-2}),
 but neither on
$\mathsf{t} \in (0,T]$ nor on $\eps>0$,} such that for every
solution
 $(  \teta, \teta_s, \uu, \chi, \xi, \eeta)$ to Problem~$\pepsmu$
  on the interval $(0,\mathsf{t})$ there holds
 \begin{align}
 \label{stimateta}
&\eps^{1/2}\|\teta\|_{L^\infty(0,\mathsf{t};V)}+
\|\teta\|_{L^2(0,\mathsf{t};V) \cap L^\infty(0,\mathsf{t};L^1(\Omega)) }
 + \| \jvstareps (\calelleveps (\teta))\|_{L^\infty (0,\mathsf{t}; L^1 (\Omega))}\leq
 C_{\mu}\,,
 \\
\label{stimatetas}
&
\begin{aligned}
\eps^{1/2}\|\teta_s\|_{L^\infty(0,\mathsf{t};\hunoc)}&+\|\teta_s\|_{L^2(0,\mathsf{t};\hunoc)
\cap L^\infty(0,\mathsf{t};L^1(\Gamma_c))} \\ & + \| \jvstarhunocps
(\calelleps (\teta_s))\|_{L^\infty (0,\mathsf{t}; L^1
(\Gamma_c))}\leq C_{\mu}\,,
\end{aligned}
\\
\label{stimachi} &\|\chi\|_{H^1(0,\mathsf{t};L^2(\Gamma_c))\cap
L^\infty(0,\mathsf{t};\hunoc)}\leq C_{\mu}\,,\\\label{stimau}
&\|\uu\|_{H^1(0,\mathsf{t};\bfw)}\leq C_{\mu}.
\end{align}
\end{lemma}
\noindent {\em Proof.} \,
 We test
 \eqref{eqIa} by
$\uu_t$:  owing to \eqref{korn_a}--\eqref{korn_b} and to
\eqref{hyp:alpha}, via the chain rule for $\widehat{\alpha}$ and the
H\"older inequality we obtain
\begin{equation}
\label{glob1}
\begin{aligned}
&\frac{C_b}{2}\int_0^t\|\uu_t\|_{\bfw}^2+\frac{C_a}{2}\|\uu(t)\|^2_{\bfw}+
\int_{\Omega} \teta \dive(\uu_t)+ \widehat{\alpha}(\uu (t))+
\int_0^t \int_{\Gamma_c} \chi\,\uu\cdot\uu_t
 \\
&\leq C \|\uu_{0}\|_W^2+\widehat{\alpha}(\uu_0) +\frac1{2 C_b}
\int_0^T \| \mathbf{F} \|_{V'}^2 \,.
\end{aligned}
\end{equation}
Next, we multiply \eqref{eqIIa} by $\chi_t$ and  integrate in time.
Again applying the chain rule to  the functional $\widehat{\beta}$,
and using~\eqref{hyp-sig}, with easy calculations we find
\begin{equation}
\label{glob4}
\begin{aligned}
\frac12\int_0^t \|\chi_t\|^2_{L^2(\Gamma_c)}&+
\frac{1}{2}\|\nabla\chi(t)\|^2_{L^2(\Gamma_c)}+\int_{\Gamma_c}\widehat{\beta}(\chi(t))
\\ &\leq C+
 \frac{1}{2}\|\nabla\chi_{0}\|^2_{L^2(\Gamma_c)}+\int_{\Gamma_c}\widehat{\beta}(\chi_0)
 -\int_0^t \int_{\Gamma_c} \lambda'(\chi)\chi_t \teta_s+
I_7 +I_8 \,,
\end{aligned}
\end{equation}
where
\begin{equation}
\label{est-intermedia} I_7= L_\sigma^2 \int_0^t \int_{\Gamma_c}
|\chi|^2 \leq 2 L_\sigma^2 T \|\chi_0 \|_{L^2 (\Gamma_c)}+ 2
L_\sigma^2 T \int_0^t \left( \int_0^s \|\chi_t \|_{L^2 (\Gamma_c)}^2
 \right)\, \dd s
\end{equation}
while,  with an easy integration by parts, we have
\begin{equation}
\label{e:inte-parts}
\begin{aligned}
I_8=-\frac12\int_0^t \int_{\Gamma_c}\chi_t |\uu|^2  & = \int_0^t
\int_{\Gamma_c} \chi \mathbf{u}_t \mathbf{u}
-\frac{1}{2}\int_{\Gamma_c} \chi(t) |\uu(t)|^2 +
\frac{1}{2}\int_{\Gamma_c} \chi_0 |\uu_0|^2\,.
\end{aligned}
\end{equation}
Furthermore, being $\widehat{\beta}$ convex,  we have
\begin{equation}
\label{beta-convex}
\begin{aligned}
\int_{\Gamma_c}\widehat{\beta}(\chi(t))&\geq -C_{1,\beta}
\|\chi(t)\|_{L^1 (\Gamma_c)} -C_{2,\beta} \\ & \geq -\eta
\|\chi(t)\|_{L^2 (\Gamma_c)}^2-C_\eta \\ & \geq -2\eta T \int_0^t \|
\chi_t\|_{L^2 (\Gamma_c)}^2 -2\eta\|\chi_0\|_{L^2 (\Gamma_c)}^2 -
C_\eta
\end{aligned}
\end{equation}
for some suitable $\eta>0$ to be specified later.
 Finally, we test
  \eqref{teta-weak-app} by
$\teta$, \eqref{teta-s-weak-app} by $\teta_s$, integrate in time and
develop in both cases the same computations as throughout
\eqref{e:chain-rule}--\eqref{coercivita}. We add the resulting
inequalities with \eqref{glob1} and \eqref{glob4}: also taking into
account~\eqref{e:inte-parts}, some terms cancel out. Furthermore,
choosing $\eta\leq 1/{8T}$ in \eqref{beta-convex} and applying the
Gronwall Lemma to deal with the integral term on the right-hand side
of \eqref{est-intermedia}, we arrive at
\begin{equation}
\label{stimaglo1}
\begin{aligned}
& \eps
\|\teta(t)\|^2_V+\|\teta\|^2_{L^2(0,t;V)}+\|\teta(t)\|_{L^1(\Omega)}
+ \eps\|\teta_s(t)\|^2_{\hunoc}+\|\teta_s\|^2_{L^2(0,t;\hunoc)}
 \\ &+\|\teta_s(t)\|_{L^1(\Gamma_c)}
+\int_0^t\int_{\Gamma_c} k(\chi)(\teta-\teta_s)^2
+\|\uu_t\|^2_{L^2(0,t;\bfw)}+\|\uu(t)\|^2_\bfw \\ & +
\int_{\Gamma_c}\chi(t)|\uu(t)|^2+\widehat\alpha(\uu(t))+\|\chi_t\|^2_{L^2(0,t;L^2(\Gamma_c))}+
\|\nabla\chi(t)\|^2_{L^2(\Gamma_c)}\\& \leq C\left(1+
\int_0^t\|h\|^2_{V'}+\int_0^t\|{\bf F}\|_{\bfw'}^2\right).
\end{aligned}
\end{equation}
for every $t \in [0,T]$. Noting that the tenth and eleventh summands
on the right-hand side are nonnegative due to \eqref{hyp:alpha} and
\eqref{A5}, we conclude. Ultimately, we also find that
\begin{equation}
\label{e:ultimate-estimates} \| \jvstareps (\calelleveps (\teta
))\|_{L^\infty (0,\mathsf{t}; L^1 (\Omega))}+ \| \jvstarhunocps
(\calelleps (\teta_s))\|_{L^\infty (0,\mathsf{t}; L^1 (\Gamma_c))}
\leq C\,,
\end{equation}
 and \eqref{stimateta}--\eqref{stimau} ensue.
 \fin
 {
 \begin{remark}
\label{e:constant-precise} \upshape As it is clear from the above
proof, constant $C_\mu$ in~\eqref{stimateta}--\eqref{stimau} depends
on $\mu$ only through $M_0^\mu$~\eqref{indip-eps-2}, i.e. the bound
for $\eps^{1/2} \| \tetazeroemu \|_{V}$ {and $\eps^{1/2} \|
\tetaessezeroemu \|_{\hunoc}$.}
 \end{remark}
 }
\subsection{Uniqueness for Problem~$\pepsmu$.}
\label{s:3.5}

We prove here the uniqueness statement in Theorem~\ref{t:3.1}. Let
us consider two families of solutions $(  \teta_i, \teta_{s,i},
\uu_i,\chi_i, \xi_i, \eeta_i)$, $i=1,2$,  to Problem~$\pepsmu$.
{Hereafter, we shall refer to the notation}
$$ \widetilde\teta=\teta_1-\teta_2,\quad
\widetilde\teta_s=\teta_{s,1}-\teta_{s,2}
,\quad\widetilde\uu=\uu_1-\uu_2, \quad\widetilde\chi=\chi_1-\chi_2
\quad\widetilde \xi=\xi_1-\xi_2, \quad\widetilde
\eeta=\eeta_1-\eeta_2.
$$
We will derive suitable contracting estimates on the solutions.
First, we subtract \eqref{teta-weak-app} written for $(\teta_2,
\teta_{s,2}, \uu_2, \chi_2)$ from \eqref{teta-weak-app} written for
$(\teta_1, \teta_{s,1}, \uu_1, \chi_1)$ and integrate in time. We
get
\begin{align}
&\eps \int_{\Omega} \widetilde{\teta}(t) \, v + \int_{\Omega}
\left(\calelleveps(\teta_1 (t))- \calelleveps(\teta_2(t))\right)v +
\eps \int_{\Omega} \nabla\widetilde{\teta}(t) \, \nabla v +
\int_{\Omega} (1 * \nabla \widetilde{\teta})(t) \, \nabla v\\
\nonumber &= \int_{\Omega} \dive \widetilde{\uu}(t) \,  v +
\int_{\Gamma_c} \Big(1
* [k(\chi_2)(\teta_2-\teta_{s,2})-k(\chi_1)(\teta_1-\teta_{s,1})]\Big) (t)\,  v
\end{align}
for all $v \in V$.
Letting $v=\widetilde{\teta}$, integrating in time and  exploiting
the monotonicity of $\calelleveps$, we obtain
\begin{align}
\label{u:1} &\eps \int_0^t \|\widetilde{\teta}\|_V^2 +\frac12 \|(1 *
\nabla \widetilde{\teta})(t)\|_{H}^2  \leq \int_0^t \int_{\Omega}
\dive \widetilde{\uu} \, \widetilde{\teta}+ \sum_{j=9}^{11}I_j.
\end{align}
In order to estimate the three latter summands, we proceed as
follows. Using H\"older's and Young's inequality (cf. also
\eqref{young}), and well-known Sobolev embeddings, we find
\begin{equation}
\label{int7}
\begin{aligned}
&I_{9}=\left|\int_0^t \int_{\Gamma_c} \Big(1
* \big[(k(\chi_2)-k(\chi_1))(\teta_1-\teta_{s,1})\big]\Big) \widetilde{\teta}\right| \\
&\leq \int_{0}^t \|\widetilde{\teta} \|_{L^2(\Gamma_c)}\, \left\|1 *
\big[(k(\chi_2)-k(\chi_1))(\teta_1-\teta_{s,1})\big]
\right\|_{L^2(\Gamma_c)}\\  &\leq C\int_{0}^t \|\widetilde{\teta}(s)
\|_{L^2(\Gamma_c)}\,
\|(k(\chi_2)-k(\chi_1))(\teta_1-\teta_{s,1})\|_{L^2(0,s;L^2(\Gamma_c))}
\, \dd s \\  &\leq C\int_{0}^t \|\widetilde{\teta}(s)
\|_{L^2(\Gamma_c)}\,
\|k(\chi_2)-k(\chi_1)\|_{L^2(0,s;L^4(\Gamma_c))}
\|\teta_1-\teta_{s,1}\|_{L^\infty(0,s;L^4(\Gamma_c))}\, \dd s\\
 & \leq \delta \int_{0}^t \|\widetilde{\teta} \|_V^2 + c_{\delta,1}
\int_{0}^t \|\widetilde{\chi} \|_{L^2(0,s;H^1(\Gamma_c))}^2\,\dd
s\,,
\end{aligned}
\end{equation}
for a suitable positive $\delta$ to be chosen later. In particular,
the positive constant $c_{\delta,1}$ also depends on the a priori
estimate~\eqref{stimateta} on
$\|\teta_1-\teta_{s,1}\|_{L^\infty(0,T;L^4(\Gamma_c))}$.
 Arguing similarly, also  taking
into account~\eqref{stimachi}  in order to estimate $\chi_2$, we
have
\begin{equation}
\label{int8}
\begin{aligned}
&I_{10}=\left|\int_0^t \int_{\Gamma_c} \Big(1
* \big[k(\chi_2)\widetilde{\teta}\big]\Big) \widetilde{\teta}\, \right| \leq \int_{0}^t \|\widetilde{\teta}
\|_{L^2(\Gamma_c)}\, \|1 * \big[ k(\chi_2)\widetilde\teta\big]
\|_{L^2(\Gamma_c)}\\  &\leq C \int_{0}^t \|\widetilde{\teta}(s)
\|_{L^2(\Gamma_c)}\, \|\chi_2+1\|_{L^\infty(0,s;L^4(\Gamma_c))}
\|\widetilde\teta\|_{L^2(0,s;L^4(\Gamma_c))}\, \dd s\\  & \leq\delta
\int_{0}^t \|\widetilde{\teta} \|_V^2 + c_{\delta,2} \int_0^t
\|\widetilde{\teta} \|_{L^2 (0,s; V)}^2 \, \dd s\,,
\end{aligned}
\end{equation}
and
\begin{equation}
\label{int9}
\begin{aligned}
&I_{11}=\left|\int_0^t \int_{\Gamma_c} \Big(1
* \big[k(\chi_2)\widetilde{\teta_s}\big]\Big) \widetilde{\teta}\,\right|
\leq \int_{0}^t \|\widetilde{\teta} \|_{L^2(\Gamma_c)}\, \|1 *
\big[k(\chi_2)\widetilde\teta_s\big] \|_{L^2(\Gamma_c)}\\
 &\leq C \int_{0}^t \|\widetilde{\teta}(s) \|_{L^2(\Gamma_c)}\,
\|\chi_2+1\|_{L^\infty(0,s;L^4(\Gamma_c))}
\|\widetilde\teta_s\|_{L^2(0,s;L^4(\Gamma_c))}\, \dd s\\  &
\leq\delta \int_{0}^t \|\widetilde{\teta} \|_V^2 + c_{\delta,3}
\int_0^t \|\widetilde{\teta_s} \|_{L^2 (0,s; H^1(\Gamma_c))}^2 \,
\dd s\,.
\end{aligned}
\end{equation}

In a similar way, we consider \eqref{teta-s-weak-app} written for
$(\teta_2, \teta_{s,2}, \chi_2)$ and for $(\teta_1, \teta_{s,1},
\chi_1)$, we take the difference,  test the resulting equation by
$\widetilde{\teta_s}$ and we integrate in time. We obtain
\begin{align}
\label{u:2} &\eps \int_0^t \|\widetilde{\teta_s}\|_{H^1(\Gamma_c)}^2
+\frac12 \|(1 * \nabla \widetilde{\teta_s})(t)\|_{L^2(\Gamma_c)}^2
\leq I_{12}+I_{13}\,,
\end{align}
where  the  terms $I_{12}$ and $I_{13}$ are estimated in the
following way. Recalling \eqref{at-most-quadratic}, we find
\begin{equation}
\label{int10}
\begin{aligned}
&I_{12}=\int_0^t \int_{\Gamma_c}
\Big(\lambda(\chi_1)-\lambda(\chi_2)\Big)\, \widetilde{\teta_s}
\\  &\leq C \int_{0}^t \|\widetilde{\chi}
\|_{L^2(\Gamma_c)}
\big(\|\chi_1\|_{L^4(\Gamma_c)}+\|\chi_2\|_{L^4(\Gamma_c)}+1\Big)\|\widetilde\teta_s
\|_{L^4(\Gamma_c)}\\  & \leq\delta \int_{0}^t
\|\widetilde{\teta_s}\|_{H^1(\Gamma_c))}^2 + c_{\delta,4} \int_0^t
\|\widetilde{\chi}\|_{L^2(\Gamma_c)}^2 \, \dd s\,.
\end{aligned}
\end{equation}
where the constant $\delta>0$ is the same as
in~\eqref{int7}--\eqref{int9}. Note that the positive constant
$c_{\delta,4}$ also depends on the estimate for
$\|\chi_1\|_{L^\infty (0,T; L^4(\Gamma_c))}$ and
$\|\chi_{2}\|_{L^\infty (0,T; L^4(\Gamma_c))}$ given
by~\eqref{stimachi}.
 Moreover, arguing as in the derivation of
\eqref{int7}--\eqref{int9}, we have
\begin{equation}
\label{int11}
\begin{aligned}
&I_{13}= \left|\int_0^t \int_{\Gamma_c} \Big(1
* \big[k(\chi_1)(\teta_1-\teta_{s,1})-k(\chi_2)(\teta_2-\teta_{s,2})
\big]\Big) \widetilde{\teta_s} \right|\\  &\leq  \delta \int_{0}^t
\|\widetilde{\teta_s} \|_{L^2(\Gamma_c)}^2 + c_{\delta,5} \int_{0}^t
\|\widetilde{\chi} \|_{L^2(0,s;H^1(\Gamma_c))}^2 \\ &  +
c_{\delta,6} \int_{0}^t \|\widetilde{\teta} \|_{L^2(0,s;V)}^2+
c_{\delta,7} \int_{0}^t \|\widetilde{\teta_s}
\|_{L^2(0,s;H^1(\Gamma_c))}^2\, \dd s\,.
\end{aligned}
\end{equation}

Now, we subtract (\ref{eqIa}),  written for
$(\teta_2,\uu_2,\chi_2,\eeta_2)$, from (\ref{eqIa}),   written for
$(\teta_1,\uu_1,\chi_1,\eeta_1)$, we test the resulting relation by
 $\widetilde\uu$ and  integrate on $(0,t)$.
 Recalling \eqref{korn_a}--\eqref{korn_b}, and the monotonicity of $\alpha$ (cf. \eqref{hyp:alpha}), we end up with
\begin{equation}
\label{u:3}
\begin{aligned}
\frac{C_b}{2}\|\widetilde\uu(t)\|^2_{W} &
+C_a\|\widetilde\uu\|^2_{L^2(0,t;W)} +\int_0^t \int_{\Gamma_c} \dive
\widetilde{\teta} \, \widetilde{\uu}
\\
&
 \leq
-\int_0^t\!\!\int_{\Gamma_c}\chi_2\,(\widetilde\uu)^2
  -\int_0^t\!\!\int_{\Gamma_c}\widetilde\chi\,\uu_1\,\widetilde\uu
\leq \int_0^t \| \widetilde\chi\|_{L^2 (\Gamma_c)} \| \uu_1 \|_{L^4
(\Gamma_c)} \| \widetilde{\uu} \|_{L^4 (\Gamma_c)}
\\
&
 \leq
C \int_0^t \| \widetilde\chi\|_{L^2 (\Gamma_c)}^2  + C \int_0^t \|
\widetilde{\bf u} \|_{W}^2 \, ,
\end{aligned}
\end{equation}
where the last two inequalities follow from  H\"older's and Young's
inequalities,
 from  Sobolev's embeddings
 and
from the fact that $\chi_2 \geq 0$ a.e. on $(0,T) \times \Gamma_c$,
due to  \eqref{A5}. In particular, the constant $C$ in \eqref{u:3}
depends on
 $\| \uu_1 \|_{L^\infty(0,T;L^4(\Gamma_c))}$ through  estimate  \eqref{stimau}.

On the other hand, let us consider the difference of \eqref{eqIIa}
written for $(\teta_{s,1},\uu_1,\chi_1,\xi_1)$ and \eqref{eqIIa} for
$(\teta_{s,2},\uu_2,\chi_2,\xi_2)$, multiply it by
$\widetilde{\chi}$ and integrate on $(0,t)\times \Gamma_c$. Taking
\eqref{A5}, \eqref{hyp-sig}, and \eqref{hyp-lambda} into account, we
get
\begin{equation}
 \label{u:4}
\begin{aligned}
&\frac{1}{2}\|\widetilde\chi(t)\|^2_{L^2(\Gamma_c)}+
 \int_0^t \| \nabla \widetilde{\chi} \|^2_{L^2(\Gamma_c)}
 \leq
 \int_0^t\!\!\int_{\Gamma_c}
(\sigma'(\chi_2)-\sigma'(\chi_1))\widetilde{\chi} \,
-\int_0^t\!\!\int_{\Gamma_c}
\lambda'(\chi_1)\,\widetilde{\teta_s}\,\widetilde{\chi}\\
& -\int_0^t\!\!\int_{\Gamma_c}
\big(\lambda'(\chi_1)-\lambda'(\chi_2)\big)\,\teta_{s,2}\widetilde{\chi}
-\frac12 \int_0^t \int_{\Gamma_c} (\uu_1 +\uu_2) \cdot
 \widetilde{\uu} \widetilde{\chi}
\\
&
 \leq L_{\sigma}\int_0^t \|\widetilde{\chi} \|^2_{L^2(\Gamma_c)}
+ C \int_{0}^t
\big(\|\chi_1\|_{L^4(\Gamma_c)}+1\big)\|\widetilde\teta_s
\|_{L^4(\Gamma_c)}\, \|\widetilde{\chi}\|_{L^2(\Gamma_c)}\\
& + C \int_{0}^t  \|\teta_{s,2}\|_{L^4(\Gamma_c)}\, \|\widetilde\chi
\|_{L^4(\Gamma_c)}\, \|\widetilde{\chi}\|_{L^2(\Gamma_c)}+
 C \int_{0}^t  \|\uu_1+\uu_2\|_{L^4(\Gamma_c)}\,\|\widetilde\uu
\|_{L^4(\Gamma_c)}\, \|\widetilde{\chi}\|_{L^2(\Gamma_c)}
\\
& \leq \delta \int_0^t \|\widetilde{\teta_s} \|^2_{H^1(\Gamma_c)}+
\delta \int_0^t \|\widetilde{\chi} \|^2_{H^1(\Gamma_c)}+ C \int_0^t
\|\widetilde{\uu} \|^2_{\WW}+ c_{\delta,8} \int_0^t
\|\widetilde{\chi} \|^2_{L^2(\Gamma_c)} \,,
\end{aligned}
\end{equation}
 where
 the constant $c_{\delta,8}$  depends on estimates
 \eqref{stimatetas}--\eqref{stimau} for the quantities
 $\| \chi_1\|_{L^{\infty}(0,T;L^4(\Gamma_c))}$, $\| \teta_{s,2}\|_{L^{\infty}(0,T;L^4(\Gamma_c))}$,
 and on $\| \uu_1 +\uu_2\|_{L^{\infty}(0,T;L^4(\Gamma_c))} $.

Finally, we add \eqref{u:1}, \eqref{u:2}, \eqref{u:3}, and
\eqref{u:4}. Noting that two terms cancel out and  taking into
account \eqref{int7}--\eqref{int9}   and
\eqref{int10}--\eqref{int11} (in which we choose,  e.g.,
$\delta=\min\{\epsi/6, 1/2\}$),
   we find
$$
\begin{aligned}
\eps \int_0^t \|\widetilde{\teta}\|_{V}^2 & +\eps \int_0^t
\|\widetilde{\teta_s}\|_{H^1(\Gamma_c)}^2 + \int_0^t
\|\widetilde{\chi}\|_{H^1(\Gamma_c)}^2
+\|\widetilde\chi(t)\|^2_{L^2(\Gamma_c)}+
C_b\|\widetilde\uu(t)\|^2_{\WW}\\
  & \leq C\Big( \int_{0}^t \|\widetilde{\teta} \|_{L^2(0,s;V)}^2+
\int_{0}^t \|\widetilde{\teta_s} \|_{L^2(0,s;H^1(\Gamma_c))}^2\\  &
+\int_{0}^t \|\widetilde{\chi} \|_{L^2(0,s;H^1(\Gamma_c))}^2\, \dd s
+ \int_0^t \|\widetilde{\chi}\|_{L^2(\Gamma_c)}^2+\int_0^t
\|\widetilde{\uu}\|_{\WW}^2\Big) \,.
\end{aligned}
$$
Thus,  by Gronwall's Lemma,
 we conclude
 that $\teta_1 =\teta_2$, $\teta_{s,1} =\teta_{s,2}$, $\uu_1 =\uu_2$, and $\chi_1=\chi_2$. A comparison in \eqref{eqIa} and \eqref{eqIIa}
 also yields  $\eeta_1 = \eeta_2$ and $\xi_1 =\xi_2 $, so that the
 uniqueness statement in Theorem~\ref{t:3.1} follows.
 \fin

\section{Proof of Theorem~\ref{t1}}
{As we mentioned in Remark~\ref{rem:explanation}, we shall pass to
the limit in Problem~$\pepsmu$ first as $\eps \searrow 0$ and
$\mu>0$ is fixed (cf. with Proposition~\ref{p:pass-1}), and then as
$\mu\searrow 0$ (see Section~\ref{s:4.2}).}
 The next result (whose  proof  is postponed to
the Appendix) concerns the construction of sequences {of initial
data $\{ \tetazeroemu\} \subset V $ and $ \{ \tetaessezeroemu\}
\subset \hunoc$ for Problem~$\pepsmu$  complying with
{\eqref{h:teta-zero}--\eqref{indip-eps-2},} and such that the
sequence of solutions to Problem~$\pepsmu$, supplemented with the
data $(\tetazeroemu,\tetaessezeroemu, \uu_0,\chi_0)$, converges to a
solution of Problem~($\mathbf{P}$) in the two consecutive limit
procedures $\eps \searrow 0$ and $\mu\searrow 0$. As, in our
construction, the data in fact depend only on the parameter $\mu>0$,
we shall denote them as $\tetazep$ and $\tetazeps $ for simplicity.}
\begin{lemma} \label{dati-iniziali}
 Assume that the initial
data $w_0$ and $z_0$ respectively comply with~\eqref{cond-teta-zero}
and~\eqref{cond-teta-esse-zero}. Then,
\begin{enumerate}
\item
 {there exists a sequence
$\{ \wzeromu \}_\mu \subset V$  fulfilling for every $\mu >0$
\begin{align}
& \label{tetazep-1-bis} \|\wzeromu\|_{H} \leq \|w_0\|_{H}\,,
\\
& \label{tetazep-2} \int_{\Omega} j^* \left(\wzeromu\right) \leq
\int_{\Omega} j^* \left(w_0\right)\,,
\end{align}
and such that
\begin{equation}
\label{e:strong-conv-wzeromu} \wzeromu \to w_0 \qquad \text{in $H$ \
\ as $\mu \searrow 0$.}
\end{equation}
 Furthermore, let us set
\begin{equation}
\label{e:precisa} \tetazep:= \calelleveps^{-1} (\wzeromu) \quad
\text{for all $\mu >0$}\,.
\end{equation}
There exists  a constant $C_{w_0}^\mu>0$, depending on $w_0$ and on
$\mu>0$ but independent of $\eps>0$, such that for all $\eps>0$
\begin{equation}
\label{tetazep-1}
 \eps^{1/2} \| \tetazep \|_{V} \leq C_{w_0}^\mu\,.
\end{equation}
}
\item
 {There exists a sequence
$\{ \zzeromu \}_\mu \subset \hunoc$  fulfilling for every $\mu >0$
\begin{align}
& \label{tetazeps-1-bis} \|\zzeromu\|_{L^2 (\Gamma_c)} \leq
\|z_0\|_{L^2(\Gamma_c)}\,,
\\
& \label{tetazeps-2} \int_{\Gamma_c} j^* \left(\zzeromu\right) \leq
\int_{\Gamma_c} j^* \left(z_0\right)\,,
\end{align}
and such that
\begin{equation}
\label{e:strong-conv-zeromu} \zzeromu \to z_0 \qquad \text{in
$L^2(\Gamma_c)$ \ \ as $\mu \searrow 0$.}
\end{equation}
Furthermore, setting
\begin{equation}
\label{e:precisas} \tetazeps:= \calelleveps^{-1} (\zzeromu) \quad
\text{for all $ \mu >0$}\,,
\end{equation}
there exists  a constant $C_{z_0}^\mu>0$, depending on $z_0$ and on
$\mu>0$ but independent of $\eps>0$, such that for all $\eps>0$
\begin{equation}
\label{tetazeps-1}
 \eps^{1/2} \| \tetazeps \|_{\hunoc} \leq C_{z_0}^\mu\,.
\end{equation}
}
\end{enumerate}
\end{lemma}
{
\begin{remark}
\label{rem:data} \upshape Note that our construction of the
sequences $\{\tetazep\}$ and $\{\tetazeps\}$ only depends on the
data $w_0$ and $z_0$. This is, ultimately, the main reason why we
have chosen to approximate the operator $\ell$ by the regularization
$\calelleveps$, instead of the usual Yosida regularization
$\ell_\mu$. Indeed, if we had used the latter approximation of
$\ell$, starting from the datum $w_0$ we should have constructed
approximate data $\tetazep \in V$  satisfying the corresponding
bound
\begin{equation}
\label{not-possible} \| (j_\mu)^* (\ell_\mu (\tetazep)) \|_{L^1
(\Omega)} \leq C \qquad \text{for all $\mu>0$,}
\end{equation}
(where $(j_\mu)^* $ is the conjugate of the Yosida approximation
$j_\mu$ of $j$),  cf. with the proof of Lemma~\ref{l:2} (clearly,
the same considerations hold for $\tetazeps$). To
deduce~\eqref{not-possible} from the condition $j^*(w_0) \in L^1
(\Omega)$,  one should virtually choose $\tetazep$ in such a way
that $\ell_\mu (\tetazep) = w_0$ (leaving aside the condition
$\tetazep \in V$). However, it is not clear to us how to carry out
this construction, since $\ell_\mu$ is not invertible. Instead,
$\calelleveps$ can be inverted, and the calculations we shall
provide in the proof of Lemma~\ref{dati-iniziali} show that the
sequence defined by~\eqref{e:precisa} complies
with~\eqref{indip-eps-1}--\eqref{indip-eps-2}.
\end{remark}
}
 {
\begin{notation}
\upshape \label{n:4} We shall denote by $$
\begin{gathered}
\label{e:cauchy-eps} \text{$\{ (\tetame,\tetasme, \ume,
\chime,\xime, \eetame) \}$ the sequence of solutions to
Problem~$\pepsmu$}
\\
\text{
 with  initial data $\{
(\tetazep,\tetazeps,\uu_0,\chi_0) \}$.}
\end{gathered}
$$
 Further,  for simplicity
we shall use the notation
\[
 \wme:=\calelleveps (\tetame), \qquad \zme:= \calelleps
(\tetasme)\,,
\]
so that, {in view of~\eqref{e:precisa}  and of~\eqref{e:precisas}
respectively,} we have
\begin{equation}
\label{e:short-hand} \wme(0)=\wzeromu, \qquad \zme(0)=\zzeromu
\qquad \text{for all $\eps>0$.}
\end{equation}
\end{notation}
}
\subsection{Passage to the limit in $\pepsmu$ as $\eps \searrow 0$}
\label{s:4.1}
\begin{proposition}
\label{p:pass-1} Let $\mu>0$ be fixed.  Under the assumptions of
Theorem~\ref{t1},
  there exists
 a (not relabeled) subsequence and  $(\tetamu,\wmu,\tetasmu, \zmu,
 \uumu,\chimu,\ximu,\eetamu)$ such that the following convergences
 hold as $\eps \searrow 0$
\begin{align}
\label{convtetae}
 &\tetame\weakto\tetamu\quad\text{in }L^2(0,T;V), \qquad
 \eps\mathcal{R}(\tetame)
 \to 0\quad\text{in }L^\infty(0,T;V'),
 \\
 &
\label{convtetase} \tetasme\weakto\tetasmu\quad\text{in }L^2(0,T;H^1
(\Gamma_c)),  \qquad \eps\mathcal{R}_{\Gamma_c}(\tetasme)
 \to 0\quad\text{in }L^\infty(0,T;\hunoc'),
\\
& \label{conveweweak}
 \wme \weak \wmu \quad \text{in }L^2(0,T;V)\,,
\\ &
\label{convewe} \eps \mathcal{R}(\tetame) + \wme \weakto \wmu \quad
\text{in }H^1(0,T;V')\,,
\\
& \label{convezeweak}
 \zme \weak \zmu \quad \text{in }L^2(0,T;\hunoc)\,,
\\ &
\label{conveze} \eps \mathcal{R}_{\Gamma_c}(\tetasme) + \zme \weakto
\zmu \quad \text{in }H^1(0,T;\hunoc')\,,
\\
& \label{convchie}
\begin{aligned} &\chime\weaksto\chimu\quad\text{in
}H^1(0,T;L^2(\Gamma_c))\cap L^\infty(0,T;\hunoc)\cap
L^2(0,T;H^2(\Gamma_c)),
\\
&\chime\rightarrow\chimu\quad\text{in
}C^0([0,T];H^{1-\delta}(\Gamma_c)) \cap
L^2(0,T;H^{2-\delta}(\Gamma_c)) \quad \text{for all $\delta>0$},
\end{aligned}
\\
 \label{convxie}
&\xime\weakto\ximu\quad\text{in }L^2(0,T;L^2(\Gamma_c)),\\
\label{convalphae}
&\eetame\weakto\eetamu\quad\text{in }L^2(0,T;\bfw'),\\
\label{convue} &
\begin{aligned}
&
\ume\weakto\uumu\quad\text{in }H^1(0,T;\bfw),\\
&\ume\to\uumu\quad\text{in }C^0([0,T];H^{1-\delta}(\Omega)^3) \quad
\text{for all $\delta>0$\,.}
\end{aligned}
\end{align}
Moreover, $\wmu$ and $\zmu$ have the further regularity
\begin{equation}
\label{e:reg-linfty}
\begin{aligned}
& \wmu \in L^\infty (0,T;H), \qquad \jvstareps\left( \wmu\right) \in
L^\infty (0,T;L^1 (\Omega))\,,
\\
&
 \zmu \in
L^\infty (0,T;L^2 (\Gamma_c)) , \qquad \jvstareps\left( \zmu\right)
\in L^\infty (0,T;L^1 (\Gamma_c))\,,
\end{aligned}
\end{equation}
and satisfy
\begin{equation}
\label{inclu-ellemu} \begin{cases} & \wmu(x,t)= \calelleveps
(\tetamu(x,t)) \quad \forae\, (x,t) \in \Omega \times (0,T), \\ &
\zmu(x,t) = \calelleveps (\tetasmu(x,t)) \quad \forae\, (x,t) \in
\Gamma_c\times (0,T)\,.
\end{cases}
\end{equation}
Further,
 the functions $(\tetamu,\wmu,\tetasmu, \zmu,
 \uumu,\chimu,\ximu,\eetamu)$ fulfil equations~\eqref{eqIa}--\eqref{bordo1} and
  \eqref{teta-weak-app}--\eqref{teta-s-weak-app} with $\eps=0$, and
the quadruple $(\wmu,\zmu,\uumu,\chimu)$ complies with the initial
conditions
\begin{equation}
\label{init-mu}
\begin{aligned}
& \wmu(0)=\wzeromu \ \ \aein \,\Omega, \quad \zmu(0)=\zzeromu \ \
\aein \ \Gamma_c, \\ & \uumu(0)=\uu_0 \ \ \aein \ \Omega, \quad
\chimu(0)=\chi_0 \ \ \aein \ \Gamma_c.
\end{aligned}
\end{equation}
\end{proposition}
\begin{notation}
\label{n:5} \upshape
 Hereafter, we shall call $\pmu$ the  boundary value
problem given by~\eqref{eqIa}--\eqref{bordo1}
and~\eqref{teta-weak-app}--\eqref{teta-s-weak-app} with $\eps=0$,
supplemented with relations~\eqref{inclu-ellemu} and the initial
conditions~\eqref{init-mu}.
\end{notation}
\noindent {\em Proof.} Let us point out that, thanks
to~\eqref{tetazep-1-bis}--\eqref{tetazep-1}
and~\eqref{tetazeps-1-bis}--\eqref{tetazeps-2}, all estimates in
Lemma~\ref{l:global-esti} hold for the sequence
$\{(\tetame,\wme,\tetasme, \zme,
 \ume,\chime,\xime,\eetame)\}_{\eps}$ with a constant $C_{\mu}$ only
 depending on the problem data,  on the initial data
 $(w_0,z_0,\uu_0,\chi_0)$, and  on $\mu>0$, but independent of
 $\eps>0$.
In order to pass to the limit in $\pepsmu$ as $\eps\searrow 0$, we
need some further estimates in addition to
\eqref{stimateta}--\eqref{stimau}. Using the latter bounds and
arguing by comparison  in \eqref{teta-weak-app} and in
\eqref{teta-s-weak-app}, we conclude that there exists $C>0$ such
that for all $\eps, \, \mu >0$
\begin{equation}\label{stimaderell}
\|\eps\mathcal{R}(\partial_t\tetame)+\partial_t \wme\|_{L^2(0,T;V')}
+ \|\eps\mathcal{R}_{\Gamma_c}(\partial_t \tetasme)+\partial_t
\zme\|_{L^2(0,T;H^1
 (\Gamma_c)')}
 \leq C\,.
\end{equation}
Similarly, a comparison in \eqref{eqIa} leads to
\begin{equation}
\label{stimaeta}
 \|\eetame\|_{L^2(0,T;\bfw')}\leq C \quad \text{for all $\eps,\, \mu >0$}.
\end{equation}
Finally, we test \eqref{eqIIa} by $\xime \in\beta(\chime)$ and get
by standard arguments
\begin{equation}
\label{stimaxi}
 \|\xime\|_{L^2(0,T;L^2(\Gamma_c))}+\|\chime\|_{L^2(0,T;H^2(\Gamma_c))}\leq
 C \quad \text{for all $\eps,\, \mu >0$}.
\end{equation}
Moreover, in view of   the Lipschitz continuity of $\calelleveps$
(cf. with~\eqref{e:lip}), joint with estimates~\eqref{stimateta}
and~\eqref{stimatetas} for $\teta$ and $\teta_s$, we conclude that
there exists a constant $\overline{C_{\mu}}>0$, depending on the
problem data and on $\mu>0$ but independent of $\eps>0$, such that
\begin{equation}
\label{e:lip-mu} \|\wme \|_{L^2 (0,T;V)} + \| \zme \|_{L^2
(0,T;\hunoc)} \leq \overline{C_{\mu}} \quad \text{for all $\eps>0$}.
\end{equation}
Combining estimates \eqref{stimateta}--\eqref{stimau} and
\eqref{stimaeta}--\eqref{e:lip-mu} with the Ascoli-Arzel\`{a}
theorem, the well-known \cite[Them.~4, Cor.~5]{Simon87}, and
standard weak compactness results,  we find that there exists an
eight-uple $(\tetamu,\wmu,\tetasmu,\zmu,\uumu,\chimu,\eetamu,\ximu)$
such that, along a suitable (not relabeled) subsequence,
convergences \eqref{convtetae}--\eqref{conveweweak},
\eqref{convezeweak},  and \eqref{convchie}--\eqref{convue} hold.
Clearly, \eqref{convewe} follows from \eqref{stimaderell} and the
second of \eqref{convtetae}. In the same way, we obtain
\eqref{conveze}. Therefore, the further
regularity~\eqref{e:reg-linfty} for $\wmu$ and $\zmu$ ensues from
the continuous embeddings $L^2 (0,T;V) \cap H^1 (0,T;V') \subset C^0
([0,T]; H)$ and $L^2 (0,T;\hunoc) \cap H^1(0,T;{\hunoc}') \subset
C^0 ([0,T]; L^2(\Gamma_c))$.} {In order to prove the second
of~\eqref{e:reg-linfty}, we exploit a \emph{Lebesgue point}
argument. Indeed,} {by~\eqref{conveweweak} and the
lower-semicontinuity of the integral functional induced by (the
convex function) $j^*_\mu$ with respect to the weak convergence in
$L^2 (0,T;V)$, we find that for all $t_0 \in (0,T)$ and $r>0$ such
that $(t_0-r,t_0+r) \subset (0,T)$
\begin{equation}
\label{e:lebesgue-point} \int_{t_0-r}^{t_0+r}\int_{\Omega}
j^*_\mu(\wmu) \leq \liminf_{\eps \searrow 0}
\left(\int_{t_0-r}^{t_0+r}\int_{\Omega} j^*_\mu (\wme) \right)\leq
2r \sup_{\eps>0}\| j^*_\mu (\wme)\|_{L^\infty (0,T;L^1(\Omega))}
\leq 2r\overline{C}\,,
\end{equation}
the latter inequality due to~\eqref{stimateta}. We now divide the
above relation by $r$ and let $r \downarrow 0$. Using  that the
Lebesgue point property  holds at almost every  $t_0 \in (0,T)$, we
 obtain the estimate
\begin{equation}
\label{e:same-constant} \| j^*_\mu (\wmu)\|_{L^\infty
(0,T;L^1(\Omega))} \leq \overline{C} \qquad \text{for all $\mu>0$,}
\end{equation}
$\overline{C} $ being the same constant as in~\eqref{stimateta}. The
analogous bound for $\jvstareps\left( \zmu\right)$ is proved in the
same way.

 Furthermore,  we
remark that
  \eqref{hyp-k}, \eqref{hyp-lambda},
\eqref{convtetae}--\eqref{convchie},  trace theorems, and Sobolev
embeddings yield that, as $\eps \searrow 0$,
\begin{equation}
\label{servissimo}
\begin{aligned}
& k(\chime)(\tetame-\tetasme)\weakto
k(\chimu)(\tetamu-\tetasmu)\quad\text{in }L^2(0,T;L^2(\Gamma_c)),
\\
& \big ( 1 * k(\chime)(\tetame-\tetasme)\big) \to \big ( 1
* k(\chimu)(\tetamu-\tetasmu)\big)\quad\text{in
}C^0([0,T];H^{-1/2}(\Gamma_c)),
\end{aligned}
\end{equation}
 as well as
\begin{equation}
\label{servissimo-2}
\begin{aligned}
&
\lambda'(\chime)\tetasme\weakto\lambda'(\chimu)\tetasmu\quad\text{in
 }L^2(0,T;L^2(\Gamma_c)),
 \\
 &
 \lambda(\chime) \weakto \lambda (\chimu) \quad \text{in } H^1
(0,T; H^1 (\Gamma_c)').
 \end{aligned}
\end{equation}

In order to pass to the limit in Problem $\pepsmu$ as $\eps \searrow
0$, we use  the above convergences and   proceed as in Section
\ref{s:3.2}.
 In particular, for \eqref{eqIa}--\eqref{bordo1} we exploit
 \eqref{convtetae}, \eqref{convchie}--\eqref{convue}
and argue in the same way as in \cite[Prop.~4.7]{bbr1}, to which we
refer the reader. Further, relying on
\eqref{convtetae}--\eqref{convchie} and the above
\eqref{servissimo}--\eqref{servissimo-2}, we  also pass to the limit
as $\eps\searrow0$ both in \eqref{teta-weak-app} and in
\eqref{teta-s-weak-app}. We  thus conclude that
$(\tetamu,\wmu,\tetasmu,\zmu,\uumu,\chimu,\eetamu,\ximu)$ satisfies
the PDE system given by~\eqref{eqIa}--\eqref{bordo1}
  \eqref{teta-weak-app}--\eqref{teta-s-weak-app} with $\eps=0$.

Finally, it remains to show~\eqref{inclu-ellemu}.
 We shall just prove the relation for $\wmu$, the
argument for $\zmu$ being completely analogous. By maximal
monotonicity of $\calelleveps$ (cf.~\cite[Lemma~1.3, p.~42]{barbu}),
it is sufficient to show that
\begin{equation}
\label{e:limsup-mu} \limsup_{\eps \searrow 0} \int_0^T \int_{\Omega}\wme \tetame
\leq \int_0^T \int_{\Omega}\wmu \tetamu\,.
\end{equation}
}which we prove integrating in time \eqref{teta-weak-app}, testing
it by $\tetame$, and again integrating on $(0,T)$. Hence, we get
$$
\begin{aligned}
&\limsup_{\eps \searrow 0}  \int_0^T\int_{\Omega} \wme \tetame  \\ &
\begin{aligned}
\leq  & \limsup_{\eps \searrow 0}\int_0^T\int_{\Omega}
\left(\eps\tetazep + w_\mu^0 -\dive(\uu_0) \right)\tetame
+\limsup_{\eps \searrow 0}\int_0^T\int_{\Omega}
\eps\nabla\tetazep \nabla\tetame
-\eps
\liminf_{\eps \searrow 0}\int_0^T  \| \tetame \|_{V}^2
\\&
  +
\limsup_{\eps \searrow 0}\int_0^T\int_{\Omega}\dive (\ume)\tetame
-\liminf_{\eps \searrow 0}\int_{\Omega} \left|\left(1* \nabla
\tetame\right)(T)\right|^2
\\&
-\liminf_{\eps \searrow 0}\int_0^T \int_{\Gamma_c} \big(1 *
k(\chime)(\tetame-\tetasme) \big) \tetame +\limsup_{\eps \searrow 0}
\int_{0}^{T} \pairing{}{}{1*h}{\tetame}
\end{aligned}
\\
&\leq
 \int_0^T\int_{\Omega}
\wmu \tetamu\,,
\end{aligned}
$$
where the second passage follows from \eqref{tetazep-1} and the
first of~~\eqref{convtetae}, from
 {the second
of~\eqref{convtetae}}, from \eqref{teta-weak} and the following
relations:
$$
\lim_{\eps \searrow 0} \int_{0}^{T}
\pairing{}{}{1*h}{\tetame}=\int_{0}^{T}
\pairing{}{}{1*h}{\tetamu},\qquad  \lim_{\eps \searrow
0}\int_0^T\int_{\Omega}\dive(\ume)\tetame=
\int_0^T\int_{\Omega}\dive(\uumu)\tetamu
$$
(due to \eqref{convtetae} and \eqref{convue}, yielding in particular
that $\dive(\ume) \to \dive(\uumu)$ in $L^\infty (0,T;
L^{6/5}(\Omega))$ as $\eps \searrow 0$),
$$
\liminf_{\eps \searrow 0}\int_{\Omega} \left|\left(1* \nabla
\tetame\right)(T)\right|^2 \geq \int_{\Omega} \left|\left(1* \nabla
\tetamu\right)(T)\right|^2
$$
(by the lower semicontinuity of the norm, since $1 * \nabla \tetame
\weaksto 1
* \nabla \tetamu$ in $L^\infty (0,T;H)$ as $\eps \searrow 0$), and, finally,
$$
\lim_{\eps \searrow 0}\int_0^T  \int_{\Gamma_c} \big(1 *
k(\chime)(\tetame-\tetasme) \big) \tetame= \int_{0}^T\int_{\Gamma_c}
\big(1 * k(\chimu)(\tetamu-\tetasmu) \big) \tetamu
$$
thanks to the second of~\eqref{servissimo} and \eqref{convtetae}
(which gives that $\tetame \weakto \tetamu$ in $L^2
(0,T;H^{1/2}(\Gamma_c))$.

{In the end,  we show~\eqref{init-mu}. The initial conditions for
$\uumu$ and $\chimu$  ensue from~\eqref{e:cauchy-eps} and
convergences~\eqref{convchie} and~\eqref{convue}. On the other hand,
thanks to~\eqref{tetazep-1} there holds
\[
\eps \mathcal{R}(\tetazep) \to 0 \qquad \text{as $\eps \searrow 0$ \
in $V'$,}
\]
hence by~\eqref{e:short-hand} and~\eqref{convewe},
$$
 \wzeromu = \lim_{\eps \searrow 0} \wme (0)=\lim_{\eps \searrow 0} \big( \eps \mathcal{R}(\tetazep) + \wme (0))
= \wmu(0)\,,
$$
where all of the above limits are meant with respect, e.g.,  to the
$H^2(\Omega)'$-topology (in fact, to the topology of any space
$\mathcal{Z}$ such that $V' \subset \mathcal{Z}$ with compact
embedding). With a completely analogous argument, we prove the
initial condition for $\zmu$ as well.
\fin  
}
\subsection{Conclusion of the proof of Theorem~\ref{t1}.}
\label{s:4.2} {We are now going to
 show that the sequence $\{(\wmu,\tetamu,\zmu,\tetasmu, \uumu,\chimu,\ximu,\eetamu)
 \}_{\mu}$ of solutions to Problem~$\pmu$ (cf. with
 Notation~\ref{n:5}) obtained in Proposition~\ref{p:pass-1} admits a
 subsequence converging as $\mu \searrow 0$ to a solution of
 Problem~$(\mathbf{P})$.
 To this aim, {we point out that there exists  a constant
 $\overline{C}>0$, \emph{independent of $\mu>0$,} such that
\begin{equation}
\label{e:stimateta-bis}
\begin{cases}
\|\tetamu\|_{L^2(0,T;V) \cap L^\infty(0,T;L^1(\Omega)) }
 + \| \jvstareps (\calelleveps (\tetamu))\|_{L^\infty (0,T; L^1 (\Omega))}\leq
 \overline{C}\,,
 \\
\|\tetasmu\|_{L^2(0,T;\hunoc) \cap L^\infty(0,T;L^1(\Gamma_c))}  +
\| \jvstarhunocps (\calelleps (\tetasmu))\|_{L^\infty (0,T; L^1
(\Gamma_c))}\leq \overline{C}\,,
\\
\|\chimu\|_{H^1(0,T;L^2(\Gamma_c))\cap L^\infty(0,T;\hunoc)}\leq
\overline{C}\,,\\
\|\uumu\|_{H^1(0,T;\bfw)}\leq \overline{C}.
\end{cases}
\end{equation}
This can be proved by testing~\eqref{teta-weak-app} (with $\eps=0$)
 by $\tetamu$,  \eqref{teta-s-weak-app} (with $\eps=0$)  by
$\tetasmu$, \eqref{eqIa} by $\uumu$, and \eqref{eqIIa} by $\chimu$,
adding the resulting equations and integrating in time. Developing
the very same calculations as for Lemma~\ref{l:global-esti}, we
conclude~\eqref{e:stimateta-bis}, whence~~\eqref{stimaxi}, as well,
for a constant independent of $\mu>0$. } Likewise, arguing by
comparison in the equations satisfied by
 $\tetamu$ and $\tetasmu$ one obtains that there exists $C>0$ such
 that for all $\mu>0$
 \begin{equation}
\label{e:stimawu} \| \partial_t \wmu\|_{L^2 (0,T;V')} + \|
\eetamu\|_{L^2 (0,T;\mathbf{W}')} + \| \partial_t \zmu\|_{L^2
(0,T;\hunoc')} \leq C\,.
 \end{equation}

Further, we are in the position of proving the following crucial
estimate
\begin{equation}
\label{e:crucial-mu} \|  \wmu\|_{L^\infty (0,T;H)} +  \|
\zmu\|_{L^\infty (0,T;L^2(\Gamma_c))} \leq C
\end{equation}
for a constant independent of $\mu>0$. Indeed, let us
test~\eqref{teta-weak-app}, with $\eps=0$, by
$\wmu=\calelleveps(\tetamu)$, \eqref{teta-s-weak-app}, with
$\eps=0$, by $\zmu=\calelleveps(\tetasmu)$, add the resulting
relations and integrate on some time interval $(0,t)$, $t \in (0,T]$
 (note that these estimates may be
performed rigorously since $\wmu \in L^2 (0,T;V)$ and $\zmu \in L^2
(0,T;\hunoc)$ for all $\mu>0$). Easy calculations lead to
\begin{equation}
\label{4.34}
\begin{aligned}
&\frac12\| \wmu(t) \|_{H}^2   + \int_0^t \int_{\Omega} \nabla
\tetamu \nabla \wmu +\frac12\| \zmu(t) \|_{L^2(\Gamma_c)}^2 +
\int_0^t \int_{\Gamma_c} \nabla \tetasmu \nabla \zmu  \\ & +
\int_0^t \int_{\Gamma_c} k(\chimu)
\left(\tetamu-\tetasmu\right)\left(\wmu-\zmu\right) \leq { \frac12
\| \wzeromu \|_{H}^2 + \frac12 \|\zzeromu \|_{L^2(\Gamma_c)}^2}+
 I_{14}+
I_{15} + I_{16}\,,
\end{aligned}
\end{equation}
where
\begin{align}
& \label{e:i-14} I_{14} = \int_0^t \int_{\Omega} |\dive(\partial_t
\uumu)| |\wmu| \leq C \| \uumu\|_{H^1 (0,T;\mathbf{W})}^2 + \frac12
\int_0^t \|\wmu\|_{H}^2\,,
\\
& \label{e:i-15}
 I_{15} = \int_0^t \int_{\Omega}|h| |\wmu|  \leq  \int_0^t \|h \|_{H}
 \|\wmu\|_{H}\,,
 \\
 &
 \label{e:i-16}
 \begin{aligned}
I_{16} &  = \int_0^t \int_{\Gamma_c}|\partial_t \lambda(\chimu)|
|\zmu| \leq C \int_{0}^t  \| \partial_t \chimu\|_{L^2(\Gamma_c)}
\left(\|\chimu\|_{L^\infty (\Gamma_c)}+1\right)
\|\zmu\|_{L^2(\Gamma_c)} \\ & \leq C \|\chimu\|_{H^1 (0,T;L^2
(\Gamma_c))}^2 + C' \int_0^t \left(\|\chimu\|_{L^\infty
(\Gamma_c)}^2+1\right) \|\zmu\|_{L^2(\Gamma_c)}^2\,,
\end{aligned}
\end{align}
the first inequality in~\eqref{e:i-16} ensuing
from~\eqref{e:allafineserve}. Now, we remark that the second and the
fourth term on the left-hand side of~\eqref{4.34} are nonnegative
thanks to~\eqref{inclu-ellemu} and the monotonicity of
$\calelleveps$. Combining the fact that $k$ takes positive values
(cf.~\eqref{hyp-k}) with the latter monotonicity argument (indeed,
it can be easily checked that the trace $ \wmu\traccia$ of $\wmu$ on
$\Gamma_c$ fulfils $ \wmu\traccia= \calelleveps(\tetamu\traccia)$),
we conclude that
 the
fifth term (in the l.h.s. of~\eqref{4.34})  is nonnegative as well.
Thus, we collect \eqref{4.34}--\eqref{e:i-16}:
recalling~\eqref{hypo-h}, estimates~\eqref{tetazep-1-bis}
and~\eqref{tetazeps-1-bis} on the initial data $\wzeromu$ and
$\zzeromu$, as well as estimates~\eqref{e:stimateta-bis}}
and~\eqref{stimaxi} (the latter yields a bound for $ \chimu$
 in $L^2 (0,T; L^\infty (\Gamma_c))$),
and applying the Gronwall lemma, we end up with~\eqref{e:crucial-mu}
as desired.

All of the above estimates,  the Ascoli-Arzel\`{a} theorem,
\cite[Them.~4, Cor.~5]{Simon87}, and standard weak compactness
results yield that there exist a  subsequence of
$\{(\wmu,\tetamu,\zmu,\tetasmu, \uumu,\chimu,\ximu,\eetamu)
 \}_{\mu}$ (which we do not relabel) and functions  $(w,\teta,z,\teta_s,\uu,\chi,\xi,\eeta)$
 for which  the first  of~\eqref{convtetae},
the first of~\eqref{convtetase},
convergences~\eqref{convchie}--\eqref{convue} and
\begin{align}
& \label{che-barba-w}
\begin{aligned}
& \wmu\weaksto w  \quad \text{in }L^\infty (0,T;H)\cap H^1
(0,T;V')\,,
\\
& \wmu\to w  \quad \text{in }C^0 ([0,T];V')\,,
\end{aligned}
\\
& \label{che-barba-z}
\begin{aligned}
 & \zmu\weaksto z  \quad \text{in }L^\infty
(0,T;L^2 (\Gamma_c))\cap H^1 (0,T;\hunoc')\,,
\\
& \zmu \to z \quad \text{in }C^0 ([0,T];\hunoc')
\end{aligned}
\end{align}
hold as $\mu \searrow 0$.  Arguing in the very same way as in the
proof of Proposition~\ref{p:pass-1}, we find that the eight-uple
$(w,\teta,z,\teta_s,\uu,\chi,\xi,\eeta)$ satisfies
equations~\eqref{teta-weak}, \eqref{teta-s-weak},
\eqref{eqIa}--\eqref{bordo1}. Furthermore,
combining~\eqref{e:strong-conv-wzeromu},
\eqref{e:strong-conv-zeromu}, \eqref{init-mu},  \eqref{che-barba-w}
and~\eqref{che-barba-z}, we conclude that the quadruple
$(w,z,\uu,\chi)$ complies with the initial
conditions~\eqref{iniw}--\eqref{iniu}.

Moreover, as already pointed out in~\cite[Sec.~4]{bfl1} (again on
the basis of~\cite[Lemma~1.3, p.~42]{barbu}), to
conclude~\eqref{inclu-w} it is sufficient to prove
\begin{equation}
\label{e:evvai} \limsup_{\mu\searrow 0} \int_0^T  \int_{\Omega}\wmu
\teta_\mu \leq \int_0^T \int_{\Omega} w\teta
\end{equation}
(analogously for proving~\eqref{inclu-z}). This can be shown by
arguing in the very same way as for~\eqref{e:limsup-mu}, so we refer
the reader to the calculations developed in the proof of
Proposition~\ref{p:pass-1}.

{Finally, we recall that estimate~\eqref{e:same-constant} holds for
a constant independent of $\mu>0$. Further, we note that, by
definition of $\calelleveps$ and convergences~\eqref{convtetae}
and~\eqref{che-barba-w}, there holds
\[
\gamma_\mu (\wmu)=\tetamu -\mu \wmu  \weakto \teta \qquad \text{in
}L^2 (0,T;H)\,,
\]
so that, by the definition~\eqref{e:max-mon-2-v} of $\gamma_\mu$,
$\wmu - \rho_\mu (\wmu) \to 0 $ in $L^2 (0,T;H)$. Hence, in view
of~\eqref{che-barba-w}
\begin{equation}
\label{e:resolvent-convergence}
 \rho_\mu (\wmu)  \weakto w\qquad
\text{in }L^2 (0,T;H)\,.
\end{equation}
Using that, by~\eqref{e:max-mon-5-v}, $j^*_\mu (\wmu) \geq j^*
(\rho_\mu (\wmu))$ a.e. in $\Omega$, we conclude that for all $t_0
\in (0,T)$ and $r>0$ such that $(t_0-r,t_0+r) \subset (0,T)$ there
holds
\[
 \int_{t_0-r}^{t_0+r} \jvstar(w) \leq
\liminf_{\mu \searrow 0} \int_{t_0-r}^{t_0+r} \jvstar(\rhoveps
(\wmu)) \leq \liminf_{\mu \searrow 0} \int_{t_0-r}^{t_0+r}
\jvstareps(\wmu) \leq 2 r \overline{C}\,,
\]
where again the first inequality follows
from~\eqref{e:resolvent-convergence} and weak lower semicontinuity
of the integral functional induced by $j^*$. With the same {Lebesgue
point} argument used  for proving~\eqref{e:same-constant}, we
infer~\eqref{reg-w-ult} (a completely analogous argument
yields~\eqref{reg-z-ult}). In the end, we point out that
\eqref{reg-w-ult} (\eqref{reg-z-ult}, respectively), \eqref{inclu-w}
(\eqref{inclu-z}, resp.), and~\eqref{hyp:gamma-2-a}
 yield an estimate for $ \teta$ in
$L^\infty (0,T;L^1 (\Omega))$ (for $ \teta_s$ in $L^\infty (0,T;L^1
(\Gamma_c))$, resp.). \fin
 }

{
\appendix
\section{Appendix} \label{sez:prel-operators}
The following result shows how the coercivity
property~\eqref{hyp:gamma-2-a} translates in terms of the  Yosida
approximation of the functional $\jvstar$.
\begin{lemma}
\label{l:-ci-serve-anticipato}
Assume~\eqref{hyp:gamma-1}--\eqref{hyp:gamma-2-a}.
 Then,
\begin{enumerate}
\item
there exists  a constant $ \overline{C}_2>0 $ ($C_1$ being the same
constant as in \eqref{hyp:gamma-2-a}) such that
\begin{equation}
\label{ci-serve-anticipata} \forall\, \mu >0, \ u \in V\, :  \ \ \mu
\| \calelleveps(u) \|_{H}^2 + \jvstareps \left(
\calelleveps(u)\right) \geq C_1 \| u \|_{L^1 (\Omega)}
-\overline{C}_2\,;
\end{equation}
\item there exists a  constant $\overline{C}_2^{*}>0 $ ($C_1$ being the same constant as in
\eqref{hyp:gamma-2-a}) such that
\begin{equation}
\label{ci-serve-anticipata-punt} \forall\, \mu >0, \ v \in  H^1
(\Gamma_c)\, : \ \  \mu \|\calelleps(v) \|_{L^2 (\Gamma_c)}^2  +
\jvstarhunocps \left(\calelleps (v) \right) \geq C_1 \|v \|_{L^1
(\Gamma_c)} - \overline{C}_2^{*}\,.
\end{equation}
\end{enumerate}
\end{lemma}
\noindent {\em Proof.} We shall just
prove~\eqref{ci-serve-anticipata}, the proof
of~\eqref{ci-serve-anticipata-punt} being completely analogous. For
a given $u \in V$, let us put  $w:=  \calelleveps(u)$: it follows
from the definition of $\calelleveps$ and from~\eqref{e:max-mon-2-v}
that $u -\mu w \in \gammav (\rhoveps (w))$, whence
\begin{equation}
\label{e:fund} \rhoveps (w) \in \ellv\left(u -\mu w \right).
\end{equation}
Therefore,  one has
\begin{equation}
\label{coerc-proof}
\begin{aligned}
 \jvstareps(w)
\geq \jvstar \left(\rhoveps (w) \right) &\geq C_1 \| u -\mu
w\|_{L^1(\Omega)} -C_2 \\ & \geq C_1 \|u \|_{L^1(\Omega)} - \mu
|\Omega|^{1/2} \|w\|_H -C_2 \\ & \geq C_1\|u \|_{L^1(\Omega)} -\mu
\| w \|_{H}^2 -\frac\mu4 C^2 -C_2\,,
\end{aligned}
\end{equation}
where the first inequality follows from~\eqref{e:max-mon-5-v}, the
second one from~\eqref{hyp:gamma-2-a}, the third one from the
H\"older inequality, and the last one from trivial computations.
Hence, we conclude~\eqref{ci-serve-anticipata}.
 \fin
\noindent
\\
We may now give the
\\
\textbf{Proof of Lemma~\ref{dati-iniziali}.} We shall just develop
the construction of the sequences $\{ \wzeromu\}$ and $\{ \tetazep
\}$, the proof of the second part of the statement being completely
analogous. For every $\mu>0$, we define $\wzeromu \in V$ as the
solution of the variational equation
\[
\int_{\Omega} \wzeromu v + \mu \int_{\Omega} \nabla\wzeromu \nabla
v= \int_{\Omega} w_0 v \qquad \text{for all $v \in V$.}
\]
Arguing in the very same way as in the proof
of~\cite[Lemma~2.4]{barbu-col-gil-gras}, we find
that~\eqref{tetazep-1-bis}--\eqref{e:strong-conv-wzeromu} hold.
By~\eqref{e:precisa}, we have \[\tetazep=\mu \wzeromu
+\gamma_{\mu}(\wzeromu),\] which ensures that $\tetazep$ is in $V$
as well, since $\gamma_\mu$ is Lipschitz continuous. Then,
\eqref{tetazep-1} immediately ensues. \fin }

                                %


                                %
                                %

\end{document}